# Forward electromagnetic induction modelling in a multilayered half-space: An open-source software tool

**Gian Piero Deidda** [1], **Patricia Díaz de Alba** [2,*], **Federica Pes** [3], and **Giuseppe Rodriguez** [4]

[1]  Department of Civil, Environmental Engineering and Architecture, University of Cagliari, 09123 Cagliari, Italy; gpdeidda@unica.it (G.P.D.)
[2]  Department of Mathematics, University of Salerno, 84084 Fisciano, Italy; pdiazdealba@unisa.it (P.D.A.)
[3]  Department of Chemistry and Industrial Chemistry, University of Pisa, 56124 Pisa, Italy; federica.pes@dcci.unipi.it (F.P.)
[4]  Department of Mathematics and Computer Science, University of Cagliari, 09124 Cagliari, Italy; rodriguez@unica.it (G.R.)
*  Correspondence: pdiazdealba@unisa.it

**Abstract:** Electromagnetic induction (EMI) techniques are widely used in geophysical surveying. Their success is mainly due to their easy and fast data acquisition, but the effectiveness of data inversion is strongly influenced by the quality of sensed data, resulting from suiting the device configuration to the physical features of the survey site. Forward modelling is an essential tool to optimize this aspect and design a successful surveying campaign. In this paper, a new software tool for forward EMI modelling is introduced. It extends and complements an existing open-source package for EMI data inversion, and includes an interactive graphical user interface. Its use is explained by a theoretical introduction and demonstrated through a simulated case study. The nonlinear data inversion issue is briefly discussed and the inversion module of the package is extended by a new regularized minimal-norm algorithm.

**Keywords:** frequency domain electromagnetic method (FDEM); electromagnetic induction (EMI); nonlinear forward modelling; nonlinear inversion; sensitivity function; MATLAB Toolbox; graphical user interface; near surface geophysics; electric conductivity; magnetic permeability.

## 1. Introduction

Electromagnetic induction (EMI) methods are proximal and remote sensing methods, among the most popular in near-surface geophysics investigation. They have been successfully used, often in combination with other geophysical techniques, in many areas spanning from environmental and hydro-geophysical investigations [1–4] to the characterization and monitoring of dismissed municipal and industrial solid waste landfills [5–8], from the quantitative evaluation of soil salinity and its spatial distribution [9–13] to soil water content monitoring [14–18], from sedimentology and soil studies [19–22] to archaeology [23–27], just to name a few.

EMI methods have been used primarily with the aim of estimating apparent electrical conductivity variability, often presented as maps, or to recover subsurface distributions of electrical conductivity, magnetic permeability [28–30], and, in some cases, the dielectric permittivity [31,32], by the inversion of the EMI responses. To these ends, it is mandatory that the physical characteristics at the survey site are such that it is possible to establish a measurable electromagnetic induction phenomenon. Since every area has its own characteristic, making it suitable or unsuitable to be investigated with a certain method, what works in some cases will not work everywhere. As Knapp and Steeples remark in [33], there are some areas where good data cannot be obtained, but there are also areas with ideal conditions for successfully investigations; for the latter it might be stated that there are areas where bad data cannot be obtained. However, the same authors warn about the risk that in areas of good data, it is also always possible to obtain bad or no data, when data acquisition parameters are not effectively designed or are not designed at all. Therefore, the results of a geophysical survey, which for the present work refers to an EMI survey, primarily rely on field data quality which, in turn, strongly depends on the quality (accuracy, resolution, and sensitivity) and appropriateness of the measuring device, as well as on the way it is used. Then, an accurate interpretation of the results, expressed as maps of apparent conductivity (and/or relative magnetic permeability) or sections of true conductivity (and/or magnetic permeability) estimated by inversion, will allow the successful achievement of the survey's goals. Forward modelling can take all these aspects into account.

As summarized in Figure 1, EMI forward modelling transforms a geological subsurface model, with its own geometry and characterized by a set of electromagnetic physical properties, into an instrumental response, which also depends on the characteristics of the measurement device (the inter-coil distance; the transmitter–receiver coil configuration; the frequency of the primary magnetic field) and on the relative position with respect to the ground it assumes during measurements (the height of the coils above the ground surface). Such modelling is essentially done after data acquisition not only to infer the properties of the ground model by inversion, as it is an indispensable part of the inverse problem (it links the device responses —i.e., data space— with the subsurface electromagnetic properties, i.e., model space; see Figure 1b), but also to facilitate the interpretation of the results, making a correlation between the observed electromagnetic responses and the expected geological models.

FORWARD MODELLING
(Direct modelling)

INVERSE MODELLING
(Backward modelling)

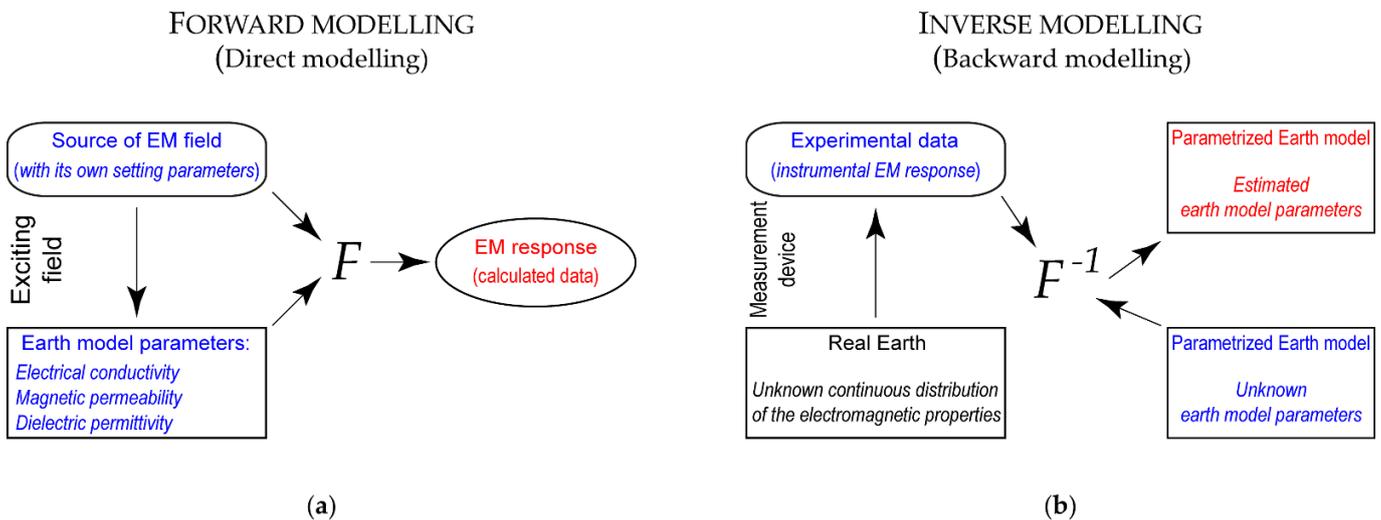

**(a)** **(b)**

**Figure 1.** (**a**) Forward modelling flowchart: $F$ is the forward modelling operator that, receiving the known subsurface material properties and the incident electromagnetic field generated by a particular instrument setting, outputs the simulated EM response (calculated data). (**b**) Inverse modelling flowchart: $F^{-1}$ is the backward modelling operator that, receiving the electromagnetic induced field measured at given positions (experimental data) and the incident electromagnetic field generated by a particular instrument setting, outputs the estimated earth model parameters.

Forward modelling should also be done before data acquisition to aid the planning of an acquisition campaign. Paraphrasing Knapp and Steeples [33], in survey and instrument design we need to start with an objective in mind, which means knowing what we wish to see. Then, we should answer the questions: what do we need to see it, and how can we get what we need to see it? That is, what characteristics (amplitude and phase) should the instrumental response have? How does it vary according to the frequency, electrical conductivity and magnetic permeability of the subsoil, the inter-coil distance, and the measurement dimension? What is the required depth of investigation? What kind of device should be used? A multi-coil instrument with different coil configurations or a multi-frequency instrument? What sensitivity should it have? Being able to run non-linear forward modelling before data acquisition would allow us to address all these issues. In addition, forward modelling is also helpful in EMI mapping for device calibration and to free measured data (apparent electric conductivity) from the "bias" introduced by the nonlinear device response, the height of the instrument, and the topography [8].

Electromagnetic induction phenomena are 3D phenomena that request full 3D forward modelling and inversion [34]. However, these require a more elaborated mathematical model of the environment and large datasets. This may be one of the reasons why 1D EMI modelling is still the most widely used system, although the literature gives examples of 3D EMI modelling and inversions [35–37].

In this work, a new Matlab-based open-source EMI 1D modelling and inversion software, FDEMtools3 (version 3.0), is introduced. FDEMtools3 comes with two graphical user interfaces (GUI), FDEMforward and FDEMinversion. The latter controls an updated version of the inversion software package described in [38], while the former drives a new sub-package devoted to EMI 1D nonlinear forward modelling, which is the focus of the present paper. Such forward modelling package has been built with the aim of providing a comprehensive tool helping to address all issues related to survey and instrument design, but also useful for an effective data inversion and a reliable data interpretation. FDEMforward is a user-friendly GUI, very well organized, and easy to access even for novice users. In addition, to make it comprehensible, as well as making the EMI method understandable to a non-specialist audience, this paper recalls the basics of electromagnetic induction and describes some mathematical aspects of the 1D forward modelling along with some key concepts of EMI methods, such as coil configurations, skin depth, induction number, sensitivity function, and depth of investigation. In this way, the present paper, together with the accompanied Matlab tool, may be viewed as a mini tutorial, ideal for teaching and training purposes. Finally, it is worth noting that the package can also be useful for advanced users since, being an open-source software, the code can be freely modified and the new functionalities can be added to meet their needs.

There are several codes that implement EMI forward and inverse modelling, some of which, such as Sim-Peg [39], FEMIC [40], and EMagPy [41], are open source, while others are commercial software, such as EM4Soil [42] and the Aarhus workbench [43]. A detailed comparison between these software codes would certainly be of great value to the Near-Surface Geophysics community, even if, at least for the moment, it is beyond the aims of the present paper. This might be done in a future work.

The structure of the present paper is as follows. Section 2 is an overview of the basic EMI theory, which Appendixes A and B complement briefly reviewing the Maxwell's equations and describing, step-by-step, the involved electromagnetic mutual induction processes. Section 3 presents the FDEMtools3 package as well as its GUIs, describing the installation process and how to use the software, by means of some numerical examples shown in Section 4. Finally, Section 5 summarizes the content of the paper.

## 2. EMI theory overview

### 2.1. Basics of electromagnetic induction

Electromagnetic induction phenomena, mostly governed by Faraday's and Ampère-Maxwell's laws, underpin the working principle of geophysical electromagnetic induction methods (Figure 1). Basically, they involve the mutual induction among three coils as shown in Figure 1b (see Appendix B for a step-by-step explanation of this mutual induction process with some mathematical details). Two of them are real metal wire coils that make up the sensors of the device. Thanks to the negligible electrical resistance of the coil-winding, from an electrical point of view they are commonly considered purely inductive circuits. Named transmitter (Tx) and receiver (Rx), they are usually placed at ground level or at a given height above it. A third further element is an imaginary coil representing a subsurface conductive magnetic body; neglecting the dielectric permittivity, it is assumed to be an RL circuit, so as to consider only resistive and inductive features of the conductive magnetic body. It is worth noting that this comes from the quasi-static approximation (see Appendix A.1), which is valid at low frequencies and when inductively induced polarization phenomena are negligible. On the contrary, at high frequencies and in the case of inductively induced polarization, an $RLC$-equivalent circuit would be more appropriate since it would also consider the capacitive features of the material in addition to the resistive and inductive ones [31,32].

An alternating current ($I_P$) passed through the loop coil Tx (Figure 2a) generates an alternating magnetic field (the primary magnetic field, $H_P$) around the loop, in-phase with the current and with the same rate of change (Figure A1c), according to Ampère-Maxwell's law. The primary magnetic field, spreading out below the ground surface, induces conduction currents and magnetization (currents of magnetization) in the conductive magnetic body. In fact, the alternating magnetic field generates a changing magnetic flux through the conductive body, which, according to Faraday's law, induces a voltage ($\mathscr{E}$) in the body, driving the so-called eddy currents ($I_{eddy}$) (Figure 2a). The induced voltage and the inductor field oscillate in sync but with a phase shift of 90° (Figure A2b,c), due to the time derivative of Faraday's law (Equation (A28)). Eddy currents may show an additional phase shift $\alpha$ (Figure A3c,d) with respect to the voltage that generates them, due to the combined effect of electrical conductivity and magnetic permeability of the target body. In fact, with reference to the coupled $LR$ equivalent circuit in Figure 2b, where resistance $R$ and inductance $L$ take into account the electrical resistivity, the magnetic permeability, and the geometry of the conductive magnetic body S, the magnitude of the additional phase lag between eddy current and voltage amounts to (Equation (A34))

$$\alpha = arctan\left(\frac{\omega L}{R}\right) = arctan(\beta), \qquad (1)$$

where $\omega$ is the radial frequency of the current in the transmitter coil (the operating radial frequency of the sensor device) and $\beta$ is the response parameter, sometimes called induction number. In circuits with only inductive loads, i.e., in perfectly conducting grounds, eddy current lags the voltage with a phase shift of 90° ($\alpha = \frac{\pi}{2}$ rad) and, therefore, with a phase lag of 180° with respect to the primary magnetic field. Otherwise, in circuits with only resistive loads, i.e., in perfectly resistive grounds, eddy current is in phase with the voltage ($\alpha = 0$ rad) and in quadrature with respect to the primary magnetic field. Due to Ampère-Maxwell's law, the time-varying eddy currents have a magnetic field associated with them, the secondary magnetic field ($H_S$) (Figure 2a), which has a phase delay of 90° plus $\alpha$ degrees with respects the primary magnetic field, as shown in Figure A4b,c. Finally, the receiver coil (Rx) simultaneously senses primary and secondary magnetic fields, measuring the voltage they induce in it according to Faraday's law. Hence, EMI devices record a complex-value electromagnetic response (Equation (A43)), usually separated into its real and imaginary parts (Equations (A44) and (A45)), which are also called In-phase (P) and Quadrature (Q) components, respectively (Figure A4c,d).

More generally, the electromagnetic response (i.e., P and Q components) is a complicated nonlinear function of many parameters, such as the electrical conductivity and the magnetic permeability of the ground, and of the technical specifications of the measuring device (the inter-coil distances, the loop-loop configurations, the height at which the sensor device operates, and the operating frequency).

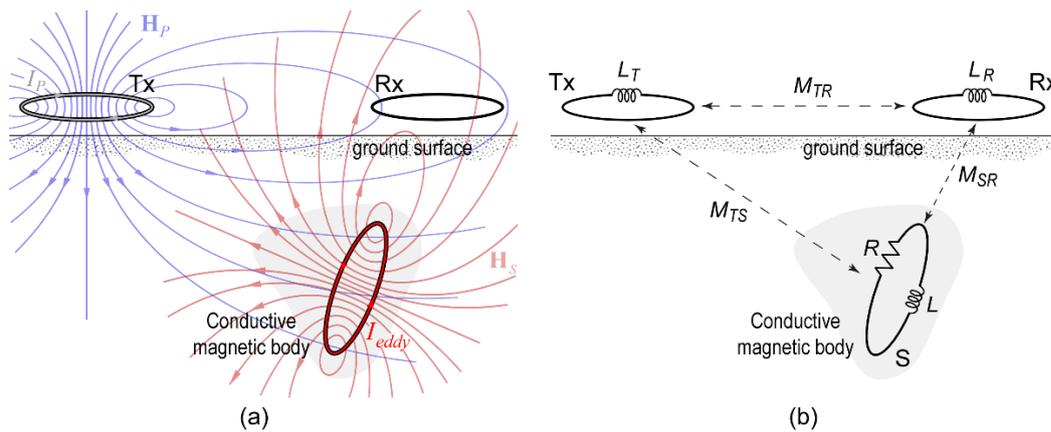

**Figure 2.** (**a**) Schematic full process of electromagnetic induction (modified from [7,44]); blue and red lines depict imaginary force lines of the primary and secondary magnetic fields, respectively. (**b**) Single-loop coupled *LR* equivalent circuits (modified from [44]); $L_T$ and $L_R$ are the self-inductance of coils Tx and Rx, respectively, while $M_{ij}$, with $i,j = T,R,S$, denotes the mutual inductance of any two of the coils.

## 2.2. 1D forward modelling

### 2.2.1. 1D layered ground model and loop-loop configurations of measuring devices

The forward modelling used to calculate the nonlinear EM response of a layered half-space for dipole source excitation is well known [45,46]. It is based on Maxwell's equations (Appendix A), suitably simplified thanks to the cylindrical symmetry of the problem. In fact, in a cylindrical coordinate system with the origin at the center of the transmitter coil, the magnetic field sensed by the receiver coil over a laterally uniform earth has no azimuthal dependence; the transmitter-receiver distance and the height above the ground surface are the only coordinates to be considered. The soil is assumed to have a layered structure with a fixed number of horizontal homogeneous layers below the subsurface, $z_1 = 0$ m (Figure 3). Each horizontal layer ranges from depth $z_k$ to $z_{k+1}$ for $k = 1, \dots, n\text{-}1$, with a thickness $d_k$, and has a conductivity $\sigma_k$ and a permeability $\mu_k$. The last layer starts at $z_n$ and is assumed to be infinite with an electrical conductivity $\sigma_n$ and magnetic permeability $\mu_n$. In the free air, above the ground surface, the conductivity is zero while the magnetic permeability is $\mu_0 = 4\pi10^{-7}$ H/m.

Modern EMI measuring devices, which are designed to collect multiple depth responses, can be grouped into multi-receiver coil systems and multi-frequency systems. The former set are endowed with multiple receiver (Rx) coils spaced at fixed distances from the transmitter (Tx) coil (Figure 3a), which usually operates at a fixed frequency; the latter set work using multiple frequencies simultaneously, usually with a fixed transmitter-receiver geometry. In addition, devices of both groups can operate at different heights above ground level, as illustrated in Figure 3a. Finally, all devices have two or more coil configurations, the most used of which are shown in Figure 3b. Table 1 lists the specifications of some commercially available devices.

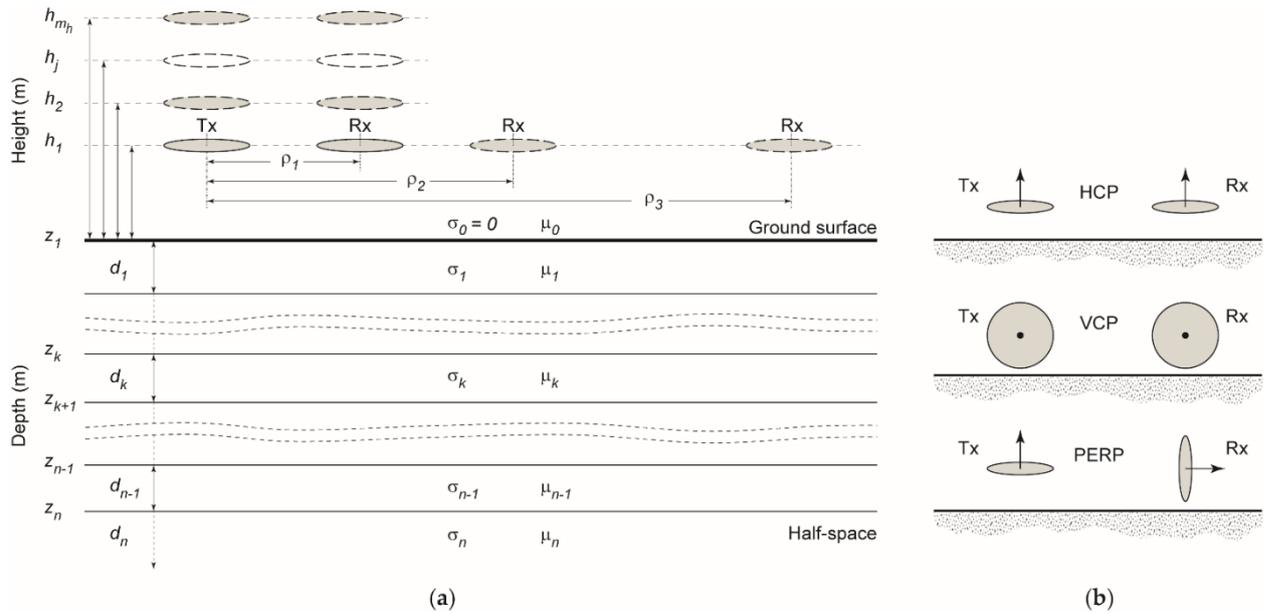

**Figure 3.** (**a**) Graphic representation discretization and parametrization for the subsoil, along with a typical measuring situation with multiple transmitter-receiver separations and/or at different heights above ground level (modified from [7]). (**b**) Common loop-loop configurations used in EMI devices: horizontal coplanar position (HCP) or vertical magnetic dipole (V); vertical coplanar position (VCP) or horizontal magnetic dipole (H); loops perpendicular to each other (PERP) or magnetic dipoles perpendicular to each other.

**Table 1.** Specifications of some commercial measuring EMI devices.

| Manufacturer | Device | Configuration | Frequency (kHz) | Coil spacing (m) | Measurement * |
|---|---|---|---|---|---|
| Gf Instruments | CMD Mini-Explorer | HCP | 30 | 0.32, 0.71, 1.18 | Q (mS/m), P (ppt) |
| | | VCP | 30 | 0.32, 0.71, 1.18 | Q (mS/m), P (ppt) |
| | CMD Explorer | HCP | 10 | 1.48, 2.82, 4.49 | Q (mS/m), P (ppt) |
| | | VCP | 10 | 1.48, 2.82, 4.49 | Q (mS/m), P (ppt) |
| | CMD DUO | HCP | 0.925 | 10, 20, 40 | Q (mS/m), P (ppt) |
| | | VCP | 0.925 | 10, 20, 40 | Q (mS/m), P (ppt) |
| Dualem Inc. | Dualem-21 | HCP | 9 | 1, 2 | Q (mS/m), P (ppt) |
| | | PERP | 9 | 1.1, 2.1 | Q (mS/m), P (ppt) |
| | Dualem-21H | HCP | 9 | 0.5, 1, 2 | Q (mS/m), P (ppt) |
| | | PERP | 9 | 0.6, 1.1, 2.1 | Q (mS/m), P (ppt) |
| | Dualem-421 | HCP | 9 | 1, 2, 4 | Q (mS/m), P (ppt) |
| | | PERP | 9 | 1.1, 2.1, 4.1 | Q (mS/m), P (ppt) |
| Geonics Limited | EM38-MK2 | HCP | 14.5 | 0.5, 1 | Q (mS/m), P (ppt) |
| | | VCP | 14.5 | 0.5, 1 | Q (mS/m), P (ppt) |
| | EM31-MK2 | HCP | 9.8 | 3.66 | Q (mS/m), P (ppt) |
| | | VCP | 9.8 | 3.66 | Q (mS/m), P (ppt) |
| Geophex Ltd. | GEM-2 | HCP | 0.03–93 | 1.66 | Q (ppm), P (ppm) |
| | | VCP | 0.03–93 | 1.66 | Q (ppm), P (ppm) |

* Q and P stand for in-Quadrature and In-Phase components, respectively; ppt and ppm are parts per thousand and parts per million, respectively.

## 2.2.2. Skin depth and induction number

An alternating current flowing in a conductor tends to distribute itself in such a way that the current density is highest near the surface of the conductor and decreases with greater depths in it. Likewise, an alternating electromagnetic field tends to concentrate near the conductor surface. In electromagnetic theory, this phenomenon is known as skin effect, the size of which is quantified by the skin depth, also called depth of penetration.

In terms of the complex wavenumber (Appendix A.1.), the skin depth in a homogeneous earth with electrical conductivity $\sigma$ and magnetic permeability $\mu$ is defined as

$$\delta = \sqrt{\frac{2}{\omega \mu \sigma}},\tag{2}$$

which represents the exponential decay of the EM-field amplitude with depth. At depth $\delta$, the EM-field amplitude has dropped by $1/e$ ($e$ is Euler's number) with respect to its value at the surface. For a $n$-layer model, the penetration in depth of the EM-fields measured at the surface ($C_1$) is solved iteratively, with a recursive formula described by the EM-response function $C_j$ [47]

$$C_j = \frac{1}{k_j} \frac{k_j C_{j+1} + tanh(k_j d_j)}{1 + k_j C_{j+1} + tanh(k_j d_j)},\tag{3}$$

where $j = n-1, n-2, \dots, 1$, $d_j$ is the thickness of the $j$th layer, $k_j = \sqrt{i\omega \mu_j \sigma_j}$ is the complex wavenumber in the $j$th layer, and $C_n = k_n$. Thus, the skin depth for a layered earth is

$$\delta = \sqrt{2|C_1|}.\tag{4}$$

This is the recursive algorithm implemented in FDEMtools3.

The induction number is another key quantity in EMI theory and practice. According to its definition in Equation (1), it depends not only on the electromagnetic properties of the conductive body but also on its geometry, even if it is difficult to evaluate unless in particular conditions. Based on the work of Grant and West [44], Ward [48] showed that the induction number depends on a linear dimension of the whole system (conductive body plus coils). In particular, for a pair of coils over a homogeneous half-space, it was shown that the induction number $\beta$ can be defined in terms of the skin depth, $\delta$, as

$$\beta = \frac{l}{\delta},\tag{5}$$

where $l$ becomes either the transmitting-receiving coil separation $\rho$ or the height $h$ of the coils above the ground, when $\rho \gg h$ or $\rho \ll h$, respectively [48].

## 2.2.3. Nonlinear forward modelling

Consider an EMI measuring device operating at angular frequency $\omega$, with coils separated by a distance $\rho$ and located at height $h$ above a 1D layered earth, like the one in Figure 3a. Assuming the configurations of the transmitting and receiving coil pair shown in Figure 3b, it would record the electromagnetic responses, defined as the ratio of the secondary ($H_S$) to the primary ($H_P$) EM field, given by

$$\begin{cases} M^{HCP}(\boldsymbol{\sigma}, \boldsymbol{\mu}; h, \omega, \rho) = -\rho^3 \int_0^\infty e^{-2\lambda h} \lambda^2 R_{\omega,0}(\lambda) J_0(\rho \lambda) d\lambda \\ M^{VCP}(\boldsymbol{\sigma}, \boldsymbol{\mu}; h, \omega, \rho) = -\rho^2 \int_0^\infty e^{-2\lambda h} \lambda R_{\omega,0}(\lambda) J_1(\rho \lambda) d\lambda \\ M^{PERP}(\boldsymbol{\sigma}, \boldsymbol{\mu}; h, \omega, \rho) = -\rho^2 \int_0^\infty e^{-2\lambda h} \lambda^2 R_{\omega,0}(\lambda) J_1(\rho \lambda) d\lambda \end{cases},\tag{6}$$

where $\boldsymbol{\sigma} = [\sigma_1, \dots, \sigma_n]^T$ and $\boldsymbol{\mu} = [\mu_1, \dots, \mu_n]^T$ represent the conductivity and the magnetic permeability vectors related to depths $d_i$, for $i = 1, \dots, n$, as shown in Figure 3, respectively, $\lambda$ is an integration variable representing the depth below the ground, normalized by the inter-coil distance $\rho$, $J_0$ and $J_1$ are Bessel functions of the first kind of zeroth and first orders, respectively, and $R_{\omega,0}(\lambda)$ is the response kernel, also called reflection factor. The reflection factor $R_{\omega,0}(\lambda)$ takes complex values and depends on the parameters that describe the layered subsurface, that is, $\sigma_k$, $\mu_k$, and $d_k$, and on the angular frequency $\omega$ and the variable of integration $\lambda$. It can be written as

$$R_{\omega,0}(\lambda) = \frac{N_0(\lambda) - Y_1(\lambda)}{N_0(\lambda) + Y_1(\lambda)},\tag{7}$$

where $Y_1(\lambda)$ and $N_0(\lambda) = \lambda/(i\omega \mu_0)$ are the surface admittance and the intrinsic admittance of the free space, respectively; in the latter, $i$ is the imaginary unit, $\omega$ is the angular frequency, and $\mu_0$ is the magnetic permeability of the free space. Setting $Y_n(\lambda) = N_n(\lambda)$, $Y_1(\lambda)$ can be obtained using Wait's back-recursive formula

$$Y_k = N_k \frac{Y_{k+1} + N_k \tan h(d_k u_k)}{N_k + Y_{k+1} \tan h(d_k u_k)}, \ k = n-1, \dots, 1, \tag{8}$$

where $d_k$ represents the $k$th layer thickness and

$$N_k = \frac{u_k(\lambda)}{i\omega\mu_k}, \tag{9}$$

is the intrinsic admittance of the $k$th layer, with

$$u_k(\lambda) = \sqrt{\lambda^2 + i\sigma_k\mu_k\omega}. \tag{10}$$

### 2.2.4. Linear approximation of the forward modelling

As a special case we recall here, for completeness, the linear case. We do this not only because it is still particularly used today in many applications but also to recall some of its limitations.

For a nonmagnetic half-space, when the coils are laid out on the ground and the operating frequency is small, the complicated relationships (6) can be well approximated in a simplified form. Under these conditions and for different coil configurations, Wait [49,50] gave a simplified expression for the secondary magnetic field as a function of the induction number $\beta$ (e.g., [50], Equations (1), (3) and (4), p. 632, for the HCP, VCP, and PERP configuration, respectively). Moreover, taking into account the relationship in [50], when the induction number is very small ($\beta \ll 1$), the imaginary part of the ratio of secondary to primary magnetic fields is linearly proportional to the half-space conductivity, $\sigma$, for both HCP and VCP coil configurations, according to [51]

$$M^{HCP} = M^{VCP} = Im\left(\frac{H_S}{H_P}\right)_{\substack{HCP \\ VCP}} = \frac{\omega\mu_0\rho^2}{4}\sigma. \tag{11}$$

This is the quadrature (Q) part of the EMI response at low induction number (LIN) condition. Most of the commercially available measuring devices incorporate the following equation to measure the apparent conductivity (as defined in [52]) directly in mS/m

$$\sigma_a = \frac{4}{\omega\mu_0\rho^2} \cdot Im\left(\frac{H_S}{H_P}\right)_{\substack{HCP \\ VCP}} \tag{12}$$

(this is why the Q component is also named LIN apparent conductivity, LIN ECa or LIN $\sigma_a$, provided that the LIN condition is met). Under these conditions, these devices also measure the in-phase part in ppt (part per thousand), which is in general very small in comparison to the quadrature (Q) component.

The general rule arising from the LIN condition can be summarized in: (1) independence between quadrature and in-phase components, i.e., the in-phase component is insignificant to generate the observations (apparent conductivity) by inversion; (2) the quadrature part is the only component directly associated to the apparent conductivity of the soil; and (3) the in-phase component is strongly connected to the magnetic susceptibility of the measured material (that is, the in-phase component is not significant for nonmagnetic materials). However, great care must be taken when using this general rule since the LIN condition only occurs when the apparent conductivity is very low (less than a few tens of mS/m) [53], which is a condition rarely met in near surface geophysics applications. As explained in Appendix B, and also pointed out in [7], P and Q components are not independent. The P component does not necessarily depend on the magnetic permeability alone, but it is mainly determined by the relative values of the inductance property with respect to the resistance property of the measured material. In fact, for a given frequency, at a fixed magnetic permeability, the P component will increase as the electrical conductivity increases, as shown in [8] (Figure S3, Supplementary Materials). Therefore, in the case of very conductive soils, the P component may be as important as the Q component and, thus, a nonlinear inversion of the complex-valued EMI response (both Q and P components, simultaneously) is needed to estimate in a good way the electrical conductivity. In addition, it is worth noting that for soils with high conductivity, the increase of the P component caused by an increase in the magnetic permeability might be hidden by the increase that P undergoes due to the electrical conductivity. This is very important when looking for magnetic targets since soil high conductivity might completely mask them, making the distinction between targeted objects and the surrounding soil a very difficult and challenging task.

### 2.3. Sensitivity function of EMI measuring devices

The sensitivity function of a measuring device is defined by the ratio between the variation of the Output and the variation of the Input, which is the quantity to be measured. For EMI devices, the sensitivity function quantifies how much the complex electromagnetic response recorded by the device is affected by a variation in the conductivity and/or permeability of a particular point (area or section) of the subsurface. The higher the absolute value of sensitivity function, the greater the influence of the subsurface region on the measurement. For a homogeneous or a layered half-space, the sensitivity, $S$, is usually calculated as a function of depth: $S = S(z)$. For each depth, the value of $S$ tells us how much measuring devices sense the changes in conductivity or magnetic permeability, given the device working parameters. In the model developed by McNeill at LIN conditions in [51], the sensitivity for all the coil orientations is a function which mainly depends on the depth, inter-coil distance, and height of the coils above the ground surface, and does not depend on the subsurface electromagnetic properties nor of the device operating frequency. These sensitivity functions are those usually provided by manufacturers in the specifications of their devices (see, for example, GF Instruments, 2020 [54]). Otherwise, when the LIN condition is not verified, the sensitivity function strongly depends on both soil conductivity and magnetic permeability, as well as on the specifications of the measuring system and its working parameters. Thus, for an EMI device with given frequency $\omega$ and inter-coil separation $\rho$, operating at height $h$ above the ground, the sensitivity function can be estimated with respect to both electric conductivity and magnetic permeability, for each of the available coil configurations, that is,

$$S_\sigma(z) = \left[\frac{\partial M^{HCP,VCP,PERP}}{\partial \sigma(z)}\right]_{h,\omega,\rho} \tag{13}$$

and

$$S_\mu(z) = \left[\frac{\partial M^{HCP,VCP,PERP}}{\partial \mu(z)}\right]_{h,\omega,\rho} \tag{14}$$

where $M$ is the complex EMI response; see Equation (6). The sensitivity function can take positive and negative values. Positive values of $S = S_{\sigma,\mu}(z)$ mean the measuring device better senses the conductive (or magnetic) materials; when $S$ takes negative values, in contrast, the device better senses poorly conductive (nonmagnetic) or resistive materials. Finally, the device no longer senses anything when the sensitivity is zero. Notice that gathering in a matrix the sensitivities of all forward responses with respect to all model parameters yields the Jacobian matrix with respect to $h$, $\omega$, and $\rho$. In this paper, such Jacobian (or sensitivity) matrix has been computed using the analytical expressions computed in [55,56], which are implemented in the FDEMtools3 package. Using this package, however, users can also optionally estimate the Jacobian through finite differences approximation.

In summary, the sensitivity function is of uppermost importance both in the survey design, as its knowledge helps to select the most appropriate and best configured measuring device, and in the solution of any nonlinear inversion, as it provides the link between the observed data and the model parameters in terms of the Jacobian matrix, allowing the update of the model vector.

## 2.4. Depth of investigation (DOI)

As described above (Section 2.1), EMI methods measure selected components of an electromagnetic field induced in a conductive soil in response to an exciting electromagnetic field generated by a device at or above the ground surface. The maximum distance (usually indicated as depth) in the subsoil within which the electromagnetic properties (electrical conductivity and magnetic permeability) of a given target in a given host produce a response that can be measured by a specific device defines the so-called Depth of Investigation (DOI). This measure plays a key role in EMI surveys as well as in other geophysical investigations. Its value is not only one of the objectives usually set in survey design, but is also crucial in the inversion processes, as it allows researchers to assess whether the inversion is data-driven or model-driven, preventing over- or misinterpretation of the inversion results [57].

Estimating the DOI is a difficult and challenging task because it depends on many variables, some of which are unknown (the real subsurface). Over the years, several estimates of it have been reported in literature. In some cases, the depth of investigation has been considered equal to the skin depth or to a multiple or a fraction of it. In other cases, it has been considered as a function of the skin depth [58–61]. Other methods, the most widespread, are based on the sensitivity function or, better, on its integral form, the cumulative sensitivity function [51,62–66]. According to these methods, the depth of investigation is the depth where its normalized integrated sensitivity function reaches a fixed threshold, such as, for example, 50, 70, 90%, or others. Without discussing the quality of these proposals, it is worth noting that all of them estimate a pseudo-depth, more or less reliable, which can anyway provide useful information. In this paper, we adopted the method described in Section 5 of [67], also based on the sensitivity function. It has been implemented in the FDEMtools3 package with the chief, though not exclusive, aim of providing the DOI as a useful output to be used in survey design before data acquisition.

## 3. Inversion algorithm

As already remarked, recent FDEM devices allow the user to record multiple simultaneous measurements with different inter-coil distances $\boldsymbol{\rho} = [\rho_1, \ldots, \rho_{m_\rho}]^T$, operating frequencies $\boldsymbol{\omega} = [\omega_1, \ldots, \omega_{m_\omega}]^T$, and heights $\boldsymbol{h} = [h_1, \ldots, h_{m_h}]^T$. In order to reconstruct the distribution of the electrical conductivity and the magnetic permeability as functions of depth from the available dataset, we denote the measurements by $b_{tij}^v$, where $t = 1, \ldots, m_\rho,\ i = 1, \ldots, m_\omega,\ j = 1, \ldots, m_h$, and $v = \{\text{HCP,VCP,PERP}\}$ represents, respectively, the vertical, horizontal, and perpendicular orientation of the coils. The data values $b_{tij}^v$ are then arranged by a suitable lexicographic ordering in a vector $\boldsymbol{b} \in \mathbb{C}^m$, where $m = \gamma m_\rho m_\omega m_h$ and $\gamma$ is the number of the orientation of the coils.

To represent the misfit between the model prediction (6) and experimental data values, we define the residual function

$$\boldsymbol{r}(\sigma, \boldsymbol{\mu}) = \boldsymbol{M}^v(\sigma, \boldsymbol{\mu}; \boldsymbol{h}, \boldsymbol{\omega}, \boldsymbol{\rho}) - \boldsymbol{b}, \tag{15}$$

where $\boldsymbol{M}^v$ is a vector containing the electromagnetic responses $M^v(\sigma, \boldsymbol{\mu}; h_j, \omega_i, \rho_t)$ for $t = 1, \ldots, m_\rho,\ i = 1, \ldots, m_\omega,\ j = 1, \ldots, m_h$. After this, we solve the following minimization problem

$$\min_{\sigma, \boldsymbol{\mu}} \frac{1}{2} \|\boldsymbol{r}(\sigma, \boldsymbol{\mu})\|_2^2, \tag{16}$$

where $\|\cdot\|_2$ denotes the Euclidean norm.

The algorithm we use for the resolution of problem (16) is based on a regularized damped Gauss–Newton method, where the regularization is achieved by a low-rank approximation of the Jacobian of the nonlinear model. Such approximation is obtained by the truncated singular value decomposition (SVD) or, depending on the adopted regularizing term, by the truncated generalized SVD (GSVD).

In recent years, this approach has been applied to the solution of (16) in various particular situations and coupled to specific techniques for evaluating the Jacobian and estimating the regularization parameter. For instance, in [55], the authors aimed at reconstructing the electrical conductivity of the soil assuming the permeability to be known considering only the quadrature component of the measurements as input, and determined the analytical expression of the Jacobian with respect to the variation of conductivity. In [68], the algorithm was adapted to devices that allow different configurations and can take simultaneous measurements. In this work, the authors also considered the possibility of processing the in-phase component of the signal.

In [56], we focused on the identification of the magnetic permeability distribution under the assumption that the conductivity was known beforehand. An important result in this work was to give the analytical expression of the Jacobian with respect to the variation of the magnetic permeability.

Later, the algorithm was updated in [67] to invert the whole complex signal sensed by the device, and to introduce a regularization term which promotes the sparsity of the solution, the so called minimum gradient support (MGS) stabilizer. The numerical algorithm was tested on real datasets collected in Sardinia (Italy), at the Molentargius Saline Regional Nature Park.

A Matlab toolbox implementing the above inversion techniques was made publicly available in [38], where it was supplemented with a graphical user interface (GUI) aiming at assisting the interested researcher in setting the parameters of the method and performing the computation. This software was used in [69] to obtain a 2D reconstruction of the electrical conductivity of a vertical section of the soil by solving a variational problem.

Besides the introduction of a tool for studying the forward modelling of the problem, the software presented in this paper slightly extends the inversion module of the package. In particular, the perpendicular orientation for the device coils has been implemented and is now available for inversion. Moreover, a new iterative algorithm based on the minimal-norm solution, presented in [70,71], has been included. It is concisely discussed in the following subsection.

### 3.1. Minimal-norm solution

In real applications, problem (16) is usually strongly underdetermined, so it does not admit a unique solution. The standard Gauss–Newton iterative algorithm, implemented in the previous version of the FDEMtools package [38], ensures unicity by imposing a regularity constraint on the iteration step, not on the solution itself. The problem of imposing a regularity constraint directly on the solution of problem (16), i.e.,

$$\begin{cases} \min_{\sigma,\mu} \|L(\sigma,\mu)\|_2^2 \\ (\sigma,\mu) \in \left\{ arg \min_{\sigma,\mu} \frac{1}{2} \|r(\sigma,\mu)\|_2^2 \right\} \end{cases} \tag{17}$$

where $L$ is a suitable regularization matrix, has been studied in [70,71].

Let us denote the solution by

$$x_k = (\sigma_k, \mu_k). \tag{18}$$

To compute the minimal-norm solution, the Gauss–Newton approximation has to be orthogonally projected, at the kth iteration, onto the null space of the Jacobian matrix $J_k = J(x_k)$.

When the regularization matrix is $L = I_{2n}$ the SVD of the matrix $J_k$ is employed. Indeed, it is well-known that the orthogonal projector may be written in terms of the SVD

$$P_{N(J_k)} = V_2 V_2^T, \tag{19}$$

where the columns of the matrix $V_2$ are orthonormal vectors spanning the null space of $J_k$. In case of $L \neq I_{2n}$ the orthogonal projector may be expressed in terms of the GSVD; see [70,71] for more details.

The resulting algorithm has been implemented in the following variants, all available in the new FDEMinversion GUI:

- MNGN

$$x_{k+1} = x_k + \alpha_k q_k - P_{N(J_k)} x_k, \tag{20}$$

  where $q_k$ is the solution of (16), and $\alpha_k$ is a step length. The damping parameter $\alpha_k$ is estimated by the Armijo-Goldstein principle. This implementation, introduced in [70], occasionally fails to converge, because the projection term may cause a considerable increase in the residual at particular iterations.

- MNGN2($\alpha$): in [71], a further damping parameter has been introduced for the projection term, through a second-order analysis of the residual $\frac{1}{2}\|r(x)\|_2^2$, as well as a strategy to automatically tune it. A simple choice is to consider a parameter $\alpha_k$ to control both terms,

$$x_{k+1} = x_k + \alpha_k(q_k - P_{N(J_k)} x_k), \tag{21}$$

  and estimate it by the Armijo–Goldstein principle.

- MNGN2($\alpha,\beta$): another possibility is to consider two independent parameters

$$x_{k+1} = x_k + \alpha_k q_k - \beta_k P_{N(J_k)} x_k. \tag{22}$$

  Additionally, in this case an automated tuning procedure has been introduced.

- MNGN2($\alpha,\beta,\delta$): this implementation is identical to the previous one, but the parameter $\beta_k$ is estimated by a different adaptive technique, which proved to be superior in the numerical simulations reported in [71].

The new implementation also allows the user to select a model profile $\bar{x}$ for the solution, which is applicable where sufficient a priori information on the physical system under investigation is available. When this does not happen, $\bar{x}$ is set to zero. The FDEMinversion GUI allows the user to select a constant profile $\bar{x}$, or to load a model from a file.

## 4. Software package

In this section, we describe the new tools available in the software package FDEMtools3, with respect to its previous version described in [38]. They consist of an extension of the forward model, the update of some of the computational routines concerning the inversion algorithm and of the corresponding graphical user interface (GUI), the introduction of a new GUI for forward modelling, and some bug corrections. In particular, the perpendicular orientation of the device coils has been integrated in the model, and a database of some of the most common commercial devices has been created. The database can be easily extended by the user by inserting the configuration of new devices, but also by introducing some non-currently available configurations, with the aim of investigating their performance.

The new Matlab toolbox FDEMtools3 is distributed as an archive file. It can be downloaded from the web page https://bugs.unica.it/cana/software (accessed on 24 March 2023). By decompressing it, a new directory "FDEMtools3" will be created. This directory must be added to the Matlab search path in order to be able to use the software from other directories. It contains the computational code as well as the user manual. The package requires the installation of P. C. Hansen's Regularization Tools package [72]; the directory of the package must be added to Matlab search path too. More detailed information on the installation process can be found in the README.txt file and in the manual.

The package contains routines for both the analysis of the forward model and the inversion procedure. Two subdirectories of the main directory, "dataforward" and "data", contain some datasets for running numerical tests with the forward and the inversion GUIs, respectively.

Table 2 lists the routines, divided in different groups, and reports a brief description for each of them. The first group "Forward Model Routines" includes the functions for computing the forward model, that is, the model prediction for a given conductivity and permeability distribution. The section "Computational Routines" contains the codes for forward and inverse modelling, including three GUIs, the "Test Scripts" are demonstration programs. Finally, the "Auxiliary Routines" list some functions needed to complete the whole process, which are unlikely to be called directly by the user, and the last group describes some further auxiliary files; see the Contents.m for details.

**Table 2.** FDEMtools3 reference.

| Forward Model Routines | |
|---|---|
| aconduct | compute the apparent conductivity |
| hratio | compute the ratio $H_S/H_P$, i.e., the device readings |
| inphase | compute the in-phase (real) component of the ratio $H_S/H_P$ |
| quadracomp | compute the quadrature (complex) component of $H_S/H_P$ |
| reflfact | compute the reflection factor |
| **Computational Routines** | |
| emsolvenlsig | reconstruct the electrical conductivity |
| emsolvenlmu | reconstruct the magnetic permeability |
| tsvdnewt | Gauss–Newton method regularized by T(G)SVD |
| jack | approximate the Jacobian matrix by finite differences |
| hankelpts | quadrature nodes for Hankel transform; see [73] |
| hankelwts | quadrature weights for Hankel transform; see [73] |
| FDEM | general graphical user interface (GUI) |
| FDEMforward | GUI for forward modelling |
| FDEMinversion | GUI for data inversion |
| **Test Scripts** | |
| drawfigures | test program for plotting the Figures in Section 5 |
| driverforward | test program for analyzing the forward model |
| driver | test program for the inversion problem |
| driver2D | test program for 2D inversion |
| **Auxiliary Routines** | |
| addnoise | add noise to data |
| chooseparam | define default parameters and test functions |
| chooseparambis | define default parameters |
| cumulativeresp | compute the cumulative response |
| fdemcomp | main code for the inversion algorithm |
| fdemdoi | compute the depth of investigation (DOI); see [67] |
| fdemplot | plot the reconstructed solution and, if available, the exact one |
| fdemprint | print information about the whole process |
| fdemsimp | compute an integral by Simpson's rule |
| forwardcomp | main code for the forward algorithm |
| mgsreg | compute the MGS regularization term; see [67] |
| morozov | choose regularization parameter by discrepancy principle |
| plotcumulative | display the cumulative response |
| plotforward | display intermediate results during forward modelling |
| plotresults | display intermediate results during inversion |
| plotsensfunc | display the sensitivity functions |
| quasihybrid | choose regularization parameter by quasi-hybrid method; see [74] |
| skindepth | compute the skin depth |
| sensitivityfunc | compute the sensitivity functions |
| **Auxiliary Files** | |
| FDEMdevices | GUI for managing the device database |
| dev.mat | data file containing the device database |
| information.pdf | file displayed by FDEMinversion |
| FDEMfwoutput.mat | data file produced by FDEMforward |
| FDEMoutput.mat | data file produced by FDEMinversion |

The most straightforward way for using the package is to run the main interface, issuing the command FDEM in the Matlab window, or running directly one of the two GUIs available: FDEMforward and FDEMinversion; see Figures 4 and 5.

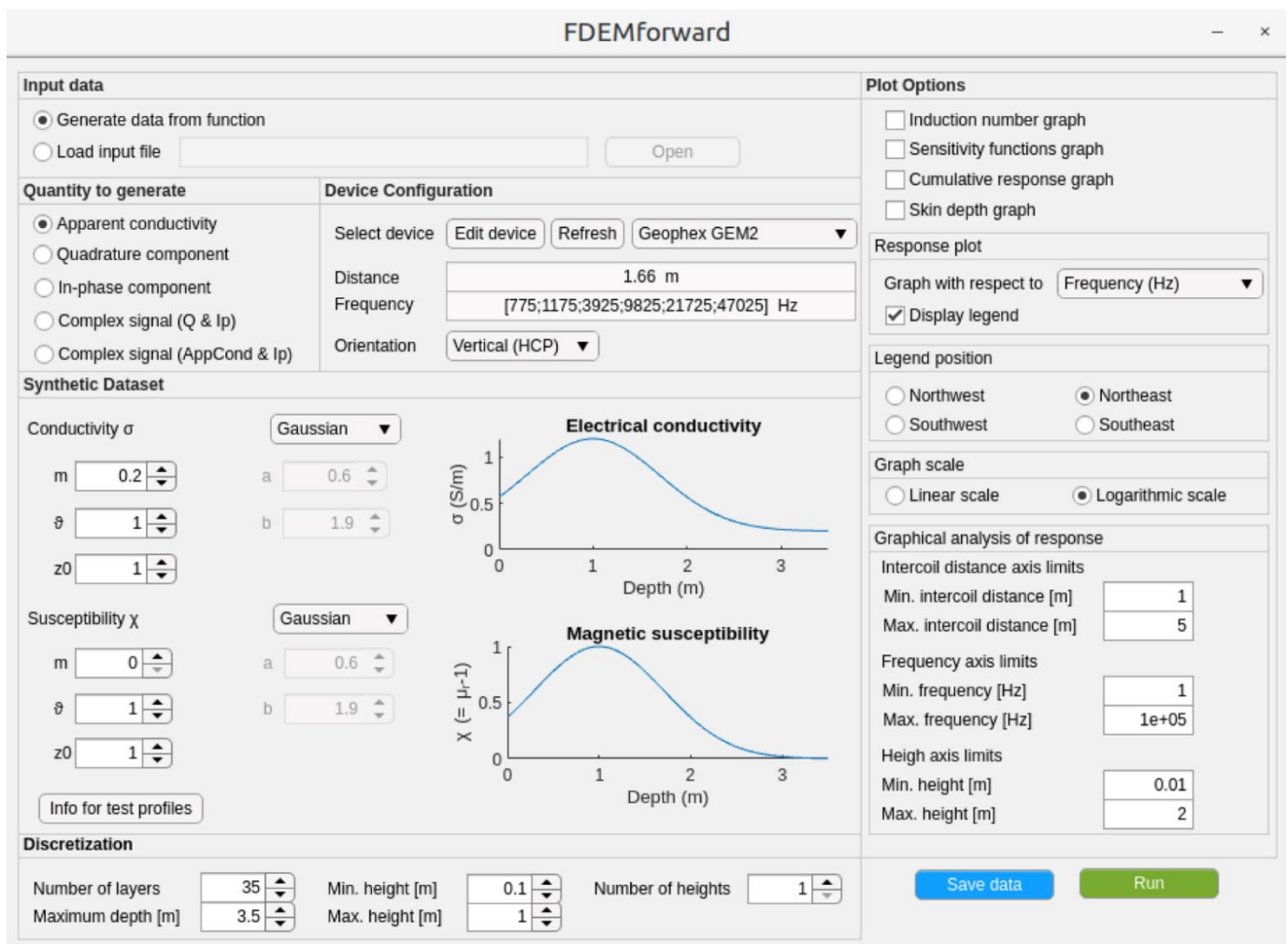

**Figure 4.** FDEMforward graphical user interface.

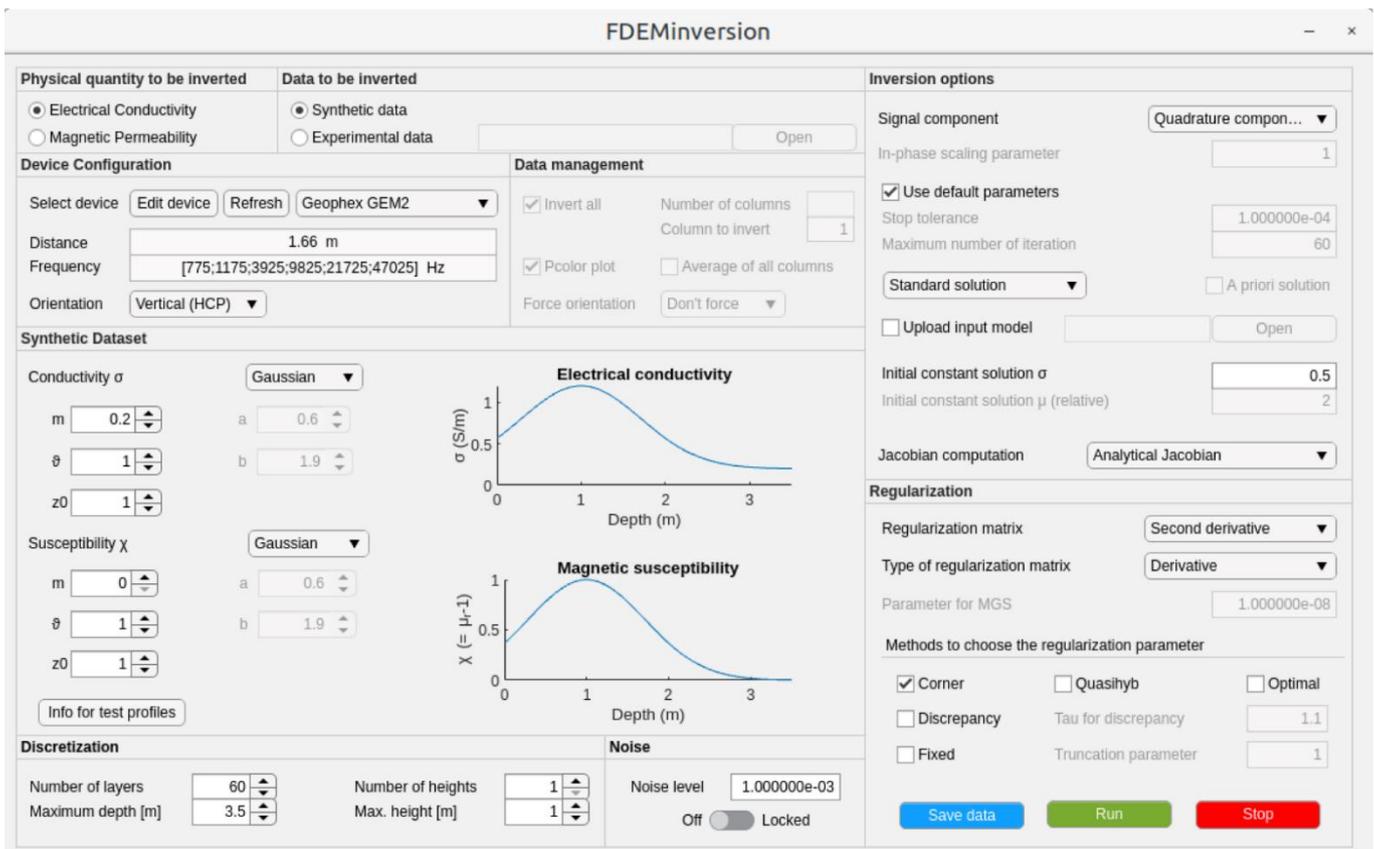

**Figure 5.** FDEMinversion graphical user interface.

Both GUIs are composed of a set of input panels that are described in detail in the user manual, and that we list here:

- FDEMforward
  o Input data;
  o Quantity to generate;
  o Device configuration;
  o Synthetic datasets;
  o Discretization
  o Plot options;
- FDEMinversion
  o Physical quantity to be inverted;
  o Data to be inverted;
  o Device configuration;
  o Data management;
  o Synthetic Dataset;
  o Discretization;
  o Noise;
  o Inversion options;
  o Regularization.

In the FDEMforward interface two buttons are available, one for running the main computation and one to save the data; in this case, a suitable data structure is used to allow the user to upload the file as an experimental dataset in FDEMinversion interface. The FDEMinversion interface contains three buttons, through which the user can start the computation, interrupt it if something goes wrong, and save the computed solution to a data file.

The computational routines can also be called directly in a Matlab script without resorting to the GUIs. This may be useful in particular situations in which, e.g., the user wants to automatize a repeated computation. We provide three example scripts for doing so: driverforward.m deals with an example of forward modelling, while driver.m and driver2D.m present two examples in which a single data column, and a set of successive data columns, are processed for inversion.

## 4. Numerical examples and discussion

This section aims to illustrate a non-exhaustive overview of some outputs of the forward modelling routines available in the FDEMtools3 package, that may be useful in survey design. For the sake of brevity, we have limited our examples to only three well-known and frequently used EMI devices, two of which, the Dualem-21H and the CMD Explorer, are multi-receiver instruments, while the other one, the GEM-2, is a multi-frequency sensor; see Table 1. Multi-receiver and multi-frequency EMI devices are able to measure the earth response at multiple depths by changing the receiver separation or frequency, respectively. Thus, for a given earth model, they can supply data suitable for resolving, by inversion, depth-related variations of electrical conductivity and/or magnetic permeability, provided that changes in receiver separation or in frequency produce in the data changes that are large enough to be measured. The Dualem-21H has one transmitter coil with a fixed frequency of 9 kHz and six receiver coils, three in a horizontal coplanar (HCP) orientation, at 0.5, 1, and 2 m from the transmitter, and three in perpendicular arrangement (PERP), at 0.6, 1.1 and 2.1 m from the transmitter. The CMD Explorer operates with one transmitter coil at a frequency of 10 kHz and has three receiver coils, spaced 1.48, 2.82, and 4.49 m from the transmitter, arranged according to the HCP or the VCP configurations. Finally, the GEM-2 contains a transmitter coil and a receiver coil separated by 1.66 m, arranged in HCP or VCP configurations, and operates in a frequency band between 30 Hz and 93 kHz, using up to ten (but usually limited to six to guarantee good signal-to-noise ratio) simultaneous frequencies.

To show and compare their responses (the signal amplitude of both the Q and P components), along with the associated sensitivities and DOIs, we have considered two three-layer earth models (Figure 6) simulating a resistive (or conductive) layer trapped between two conductive (or resistive) ones, representing targets typically found in environmental, engineering and archaeological investigations, such as contaminant plumes, foundations, archaeological structures (e.g., walls, stone built remains, ditches, tombs, and so on). In detail, the 1D earth model consists of a top layer of nonmagnetic material, with a fixed conductivity of 0.1 S/m, an intermediate layer having a relative magnetic permeability of 1.01 with a conductivity between 0.001 S/m (low conductivity case) and 2 S/m (high conductivity case), and a third layer (a half-space) with a conductivity of 0.01 S/m and a relative magnetic permeability of 1.005. The thickness of the first layer is 1.5 m while that of the middle layer is 1 m. The magnetic permeability of the intermediate layer is probably a little higher than that usually found in real soils, but it has been used to better highlight the effects that magnetic materials might have on EMI responses. In the following, the earth models with the least conductive and most conductive middle layer will be named M1 (Figure 6a) and M2 (Figure 6b), respectively.

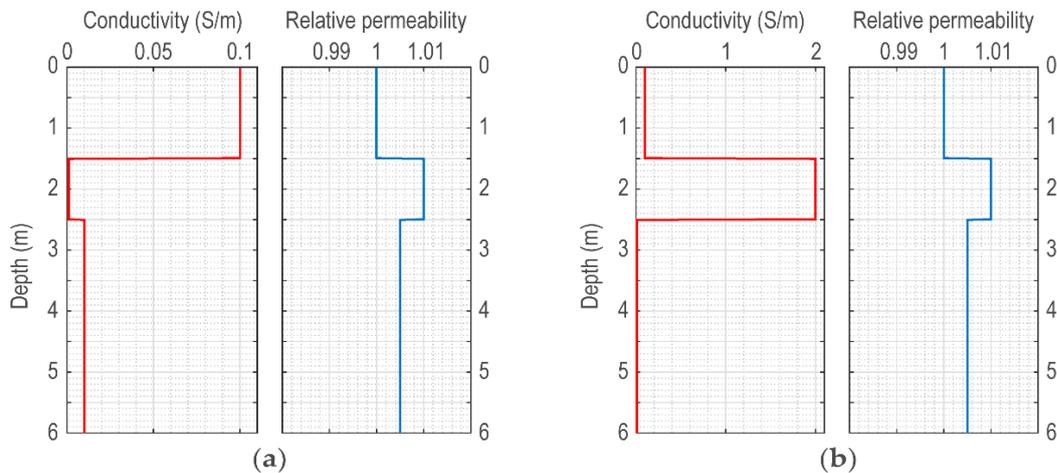

**Figure 6.** (**a**) Earth model M1; (**b**) earth model M2. Models differ only in the electrical conductivity of the middle layer.

Table 3 lists the skin depth and the induction number values arising from the combination of the device characteristics with the physical properties of the earth models. Such values have been estimated using Equations (4) and (5) and represent some optional outputs that the user can get running the FDEMforward tool.

**Table 3.** Skin depths and induction numbers.

| Device | $\rho$ (m) | $f$ (Hz) | $\delta_1$ (m) | $\beta_1$ | $\delta_2$ (m) | $\beta_2$ |
|---|---|---|---|---|---|---|
| Dualem-21H [1] | 0.5 (0.6) | 9,000 | 41.4 | 0.012 (0.015) | 8.8 | 0.057 (0.068) |
| | 1 (1.1) | 9,000 | 41.4 | 0.024 (0.027) | 8.8 | 0.113 (0.124) |
| | 2 (2.1) | 9,000 | 41.4 | 0.048 (0.051) | 8.8 | 0.226 (0.237) |
| CMD Explorer | 1.48 | 10,000 | 38.8 | 0.038 | 8.1 | 0.183 |
| | 2.82 | 10,000 | 38.8 | 0.073 | 8.1 | 0.348 |
| | 4.49 | 10,000 | 38.8 | 0.116 | 8.1 | 0.554 |
| GEM-2 | 1.66 | 1,275 | 127.9 | 0.013 | 48.2 | 0.034 |
| | 1.66 | 4,250 | 64.9 | 0.026 | 16.9 | 0.098 |
| | 1.66 | 12,525 | 33.7 | 0.049 | 6.8 | 0.246 |
| | 1.66 | 28,725 | 19.6 | 0.085 | 3.9 | 0.427 |
| | 1.66 | 54,150 | 12.6 | 0.132 | 3.1 | 0.544 |
| | 1.66 | 82,150 | 9.3 | 0.179 | 2.8 | 0.592 |

[1] The values in parentheses are for the PERP configuration. $\rho$ is the inter-coil distance and $f$ the operating frequency of each device; $\delta$ and $\beta$ are the skin depth and the induction number, respectively. Subscripts 1 and 2 refer to earth models M1 and M2.

EMI sensors can be hand-carried by a person using shoulder-harnesses or harness straps in station-by-station on-ground measurements or in a continuous-recording walking survey, but they can also be mounted on a sled or cart to be towed by a small all-terrain vehicle or tractor. However, it is worth noting that changing the way a sensor is used, whose choice is usually dictated by the desired speed of investigation, changes its operating height, which is a survey parameter that should be carefully selected as it may be a decisive factor for the success of the survey, for both imaging and mapping purposes. In fact, varying the probe height changes the depth of penetration of EMI devices, so that measurements investigate different and overlapping soil volumes [63]. This is the reason why, even when a device with a single frequency and a single receiver is used, data recorded at different instrumental heights can be inverted to get quantitative estimates of depth variations in true electrical conductivity [51,55,75–79]; the greater the effect of the height, the better the inverted result will be. The effect of the operating height remains important also for multiple depth responses collected with a multi-receiver device. To recover good estimates of conductivity with depth by inverting data measured at multiple inter-coil spacings, the device should operate at such a height that the values recorded by each coil are well separated. Concerning EMI mapping surveys, on the other hand, it is worth noting that an increase in the probe height usually lowers the amplitude of the measured response, causing the drawbacks discussed in Deidda et al. (2022) [8].

Therefore, knowing a priori how EMI responses vary as the operating height of the sensor changes, as shown by the graphs in Figures 7–10, may be very useful in survey design. For example, looking at the response of the Dualem-21H above the M1 model (Figure 7), an operating height of 0.9 m (the height the sensor would have by carrying it with a harness strap) would provide well-separated quadrature values for both HCP and PERP configurations, well suited to be inverted. This is not the case for the responses (Figure 8) the Explorer would have recorded when operating at 0.9 m above the earth model M2. In fact, the HCP quadrature values for the inter-coil distances of 1.48 m and 2.82 m (Figure 8a), as well as the VCP quadrature values for the inter-coil distances of 2.82 m and 4.49 m (Figure 8b), differ by less than 2 mS/m, which in practice may be a value smaller than the noise level. Thus, with reference to earth model M1, it turns out that the Dualem-21H would operate better at heights greater than about 0.8 m, while the Explorer would provide good data when operating directly on the ground surface. Inspecting the responses above model M2 (Figures 9 and 10), it appears that both devices would record very good data at all the considered operating heights, except those from 0.3 m to 0.5 m and from 1 m to 1.4 m for the Explorer HCP quadrature response (Figure 10a).

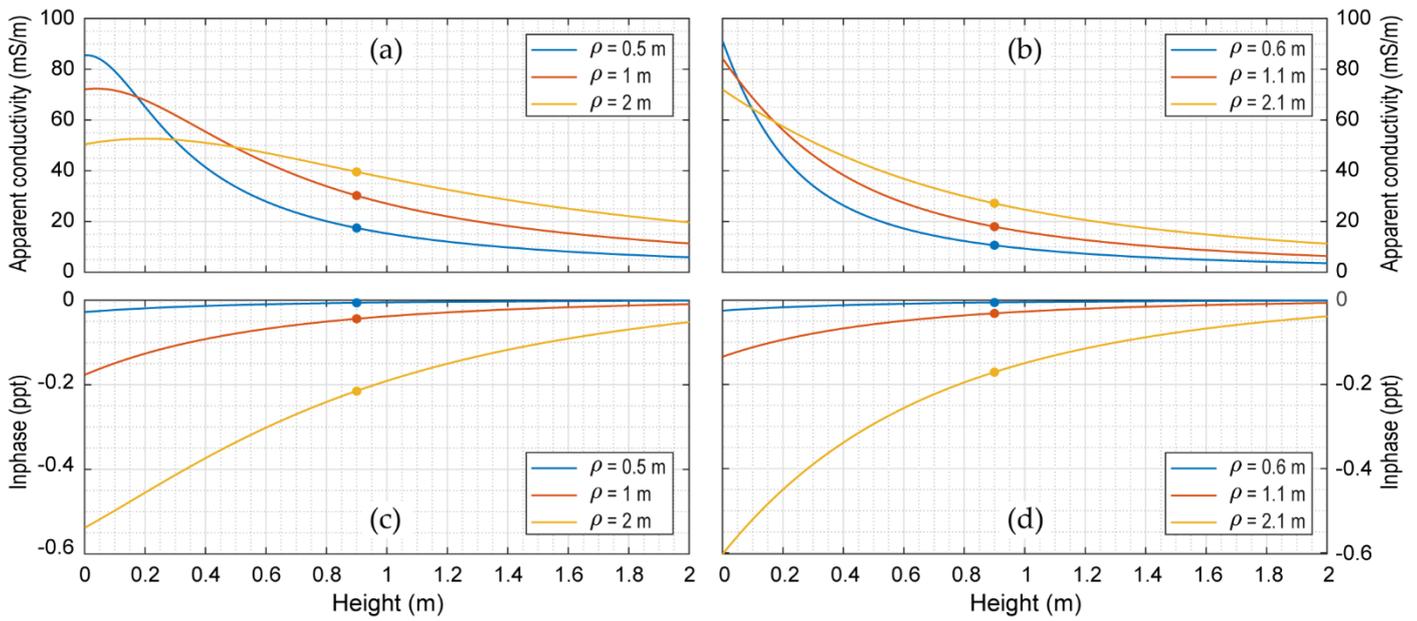

**Figure 7.** Electromagnetic response of the Dualem21H above the earth model M1. (**a**,**b**) are the simulated HCP and PERP quadrature (Q) responses, respectively, both expressed as apparent conductivity in mS/m; (**c**,**d**) are the HCP and PERP in-phase (P) responses, respectively. Dots indicate the response values at the probe height of 0.9 m, which is a frequently used operating height for both devices.

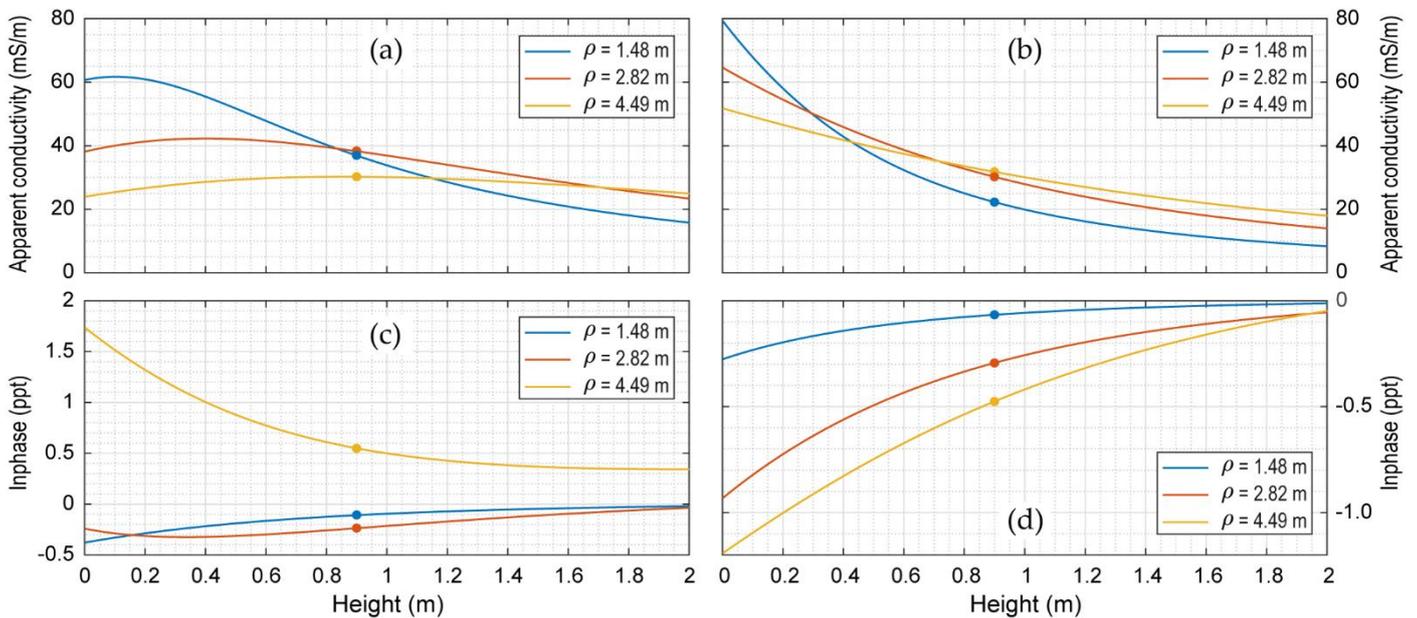

**Figure 8.** Electromagnetic response of the Explorer above the earth model M1. (**a**,**b**) are the simulated HCP and VCP quadrature (Q) responses, respectively, both expressed as apparent conductivity in mS/m; (**c**,**d**) are the HCP and VCP in-phase (P) responses, respectively. Dots indicate the response values at the probe height of 0.9 m, which is a frequently used operating height for both devices.

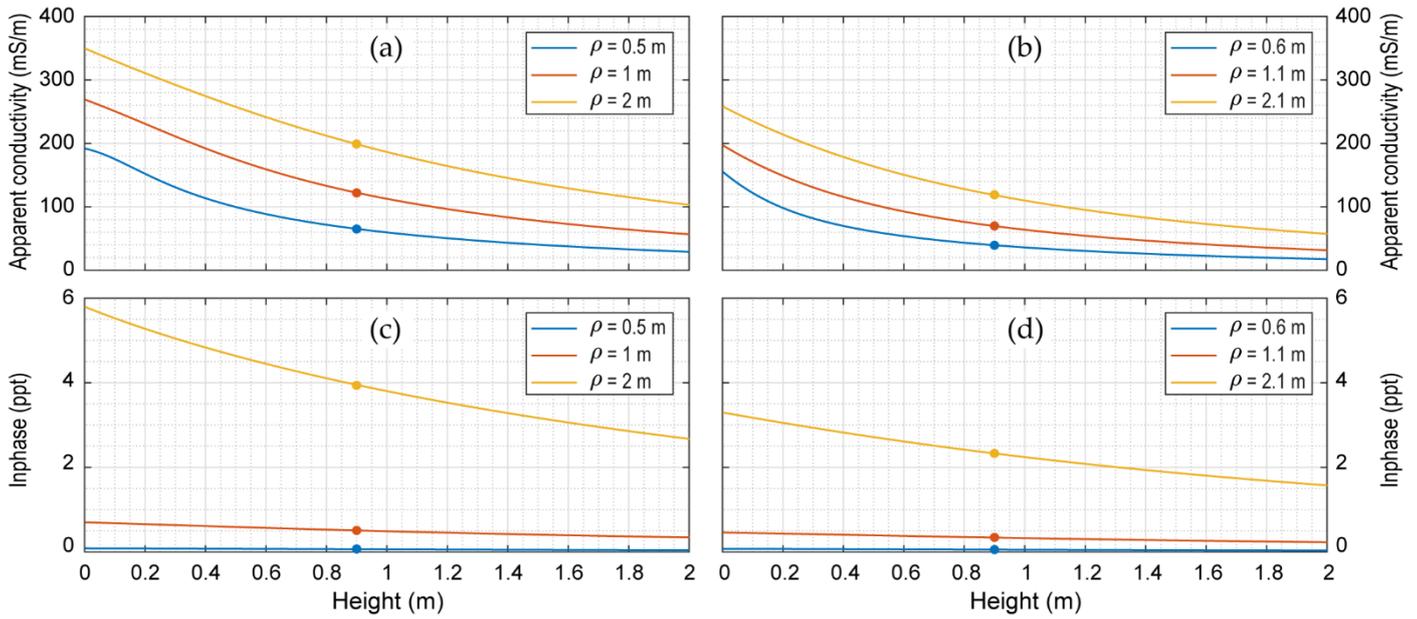

**Figure 9.** Electromagnetic response of the Dualem21H above the earth model M2. (**a**,**b**) are the simulated HCP and PERP quadrature (Q) responses, respectively, both expressed as apparent conductivity in mS/m; (**c**,**d**) are the HCP and PERP in-phase (P) responses, respectively. Dots indicate the response values at the probe height of 0.9 m, which is a frequently used operating height for both devices.

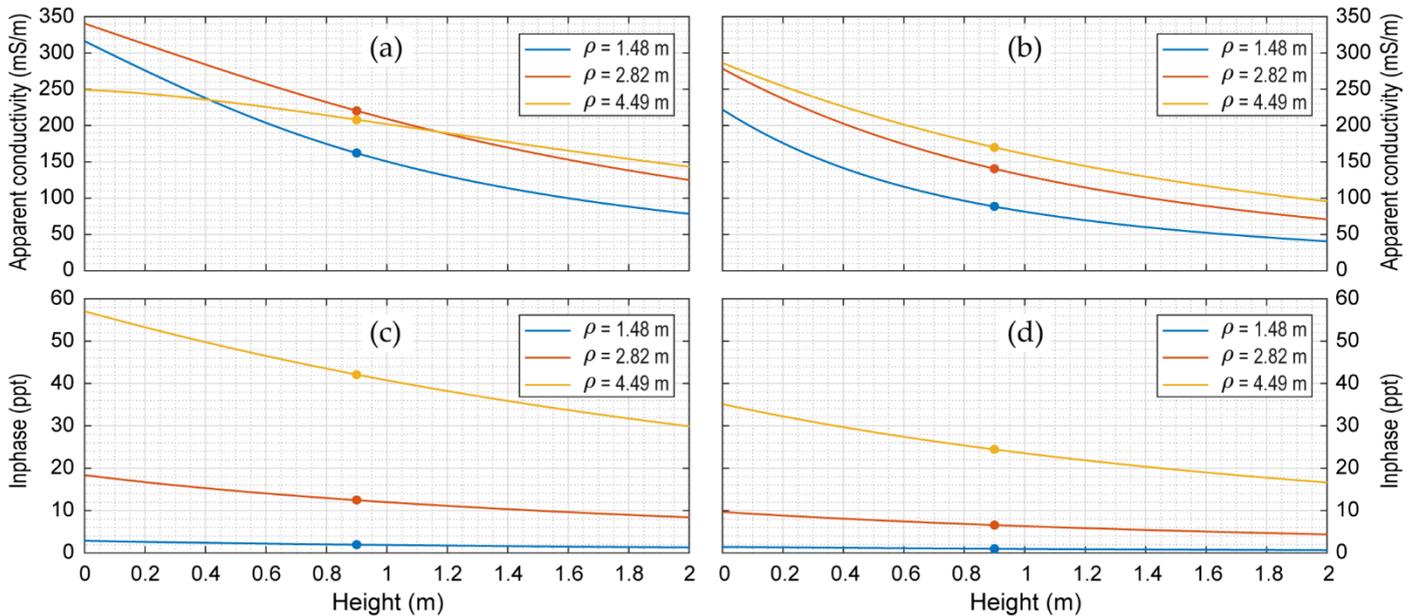

**Figure 10.** Electromagnetic response of the Explorer above the earth model M2. (**a**,**b**) are the simulated HCP and VCP quadrature (Q) responses, respectively, both expressed as apparent conductivity in mS/m; (**c**,**d**) are the HCP and VCP in-phase (P) responses, respectively. Dots indicate the response values at the probe height of 0.9 m, which is a frequently used operating height for both devices.

Inspecting the in-phase component of the responses over the earth model M1, Figures 7c,d and 8c,d show that for both multi-receiver devices the values are always very small and negative, with the only exception being the in-phase component of the Explorer at 4.49-m in the HCP configuration, whose values are positive. The presence of negative values of the in-phase component is definitely linked to the presence of susceptible materials. In fact, by running the forward modelling over the earth model M1 with relative magnetic permeability equal to 1 for all layers, the values of the in-phase component become positive. This suggests that negative values in the in-phase component may indicate the presence of susceptible materials, at least when the signal amplitude is sufficiently large to be clearly above the noise level. On some occasions, large values of the in-phase component are deemed to be evidence of the presence of magnetic materials. As already observed in Section 2.1, this is not always the case. For example, the Explorer's responses over the earth model M2 show an in-phase component (Figures 9c,d and 10c,d) with larger values than for the model M1, reaching up to about 56 ppt (Figure 10c). As the two earth models differ only for their electrical conductivity, the strong increase of the in-phase component values is due only to the electrical conductivity, and not to the magnetic permeability.

Figure 11 presents the complex electromagnetic response of the multi-frequency GEM-2 system over earth models M1 (Figure 11a) and M2 (Figure 11b). Both Q and P components of the response function are shown as a function of frequency, in the range of 30 Hz to 93 kHz. In addition, to show how the operating height affects the response values, the response has been estimated at the heights of 0.2 and 0.9 m, which are the usual heights that the device would have when hand-carried with a shoulder-strap or harness.

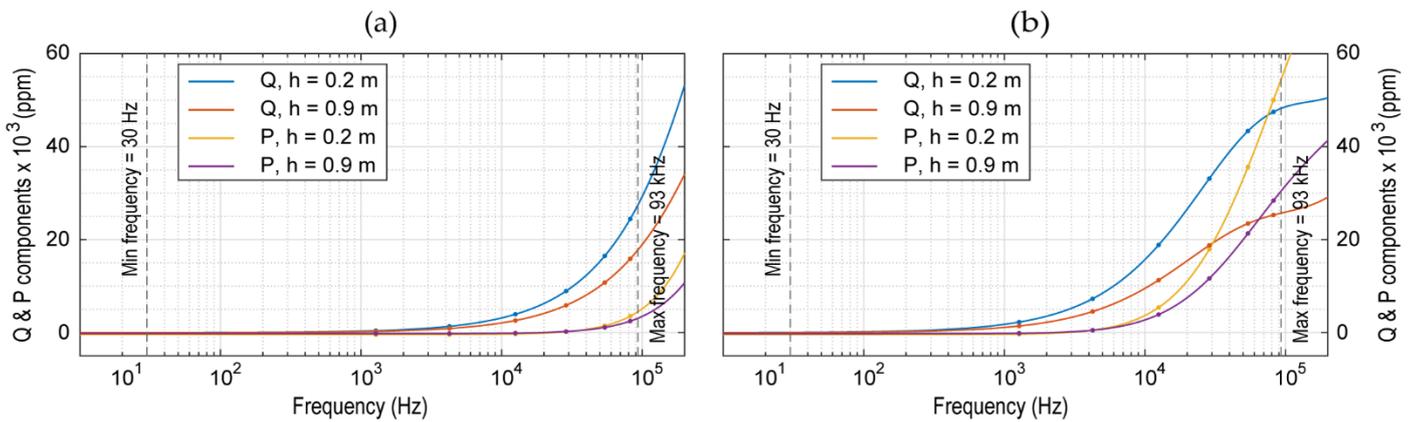

**Figure 11.** Electromagnetic response of GEM-2. (**a**) Quadrature and in-phase components of responses at heights of 0.2 and 0.9 m above ground surface of earth model M1; (**b**) Quadrature and in-phase components of responses at heights of 0.2 and 0.9 m above ground surface of earth model M2. Dots indicate the response values for a set of six selectable operating frequencies among those currently available for the GEM-2 (minimum frequency = 30 Hz; maximum frequency = 93 kHz).

As Figure 11 shows, for both earth models M1 and M2, the signal amplitude of Q and P components, very low at low frequencies, increases as the frequency increases. In addition, it is very clear that this increase is sharper over the more conductive earth model M2 and, for both models M1 or M2, at small operating heights. The small responses at low frequencies are related to the corresponding low induction numbers, defined as $\beta \leq 0.02$ in [80] (Figure 12). As explained in Appendix B (Figure A5), this means that both complex response functions become purely imaginary (resistive limit) as the induction numbers approach zero. In other words, this also means that inductive phenomena are negligible at low frequencies (low induction numbers), resulting in a small EMI response and a marginal frequency dependence, which render data inversion unfeasible. Therefore, to obtain the most useful information about earth models M1 and M2, obtaining responses that are strong, frequency dependent, and suitable for data inversion, the GEM-2 should be configured with the widest possible set of frequencies [81] to operate over a range of moderate induction number (defined as $0.02 < \beta < 1$) (Figure 12). A possible set of frequencies meeting these requirements is listed in Table 3 and shown in Figures 11 and 12. However, it is worth noting that when using this set of frequencies, despite the fact that both responses (at 0.2 m and 0.9 m) over the earth model M1 are frequency dependent, only the one estimated at 0.2 m still has acceptable signal amplitudes.

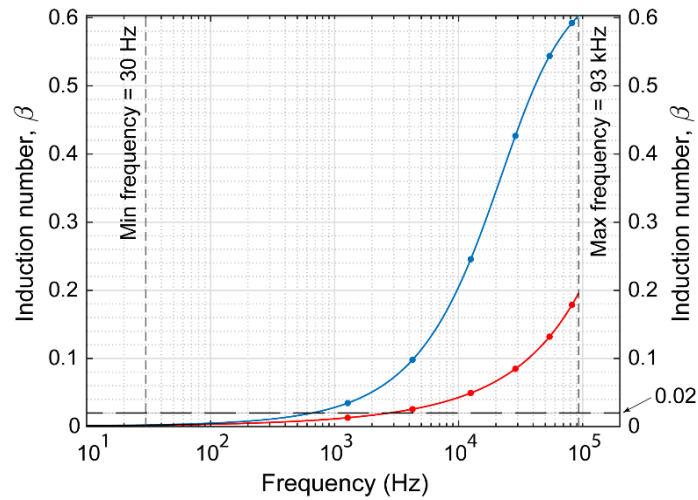

**Figure 12.** Induction numbers spanned over earth models M1 (red curve) and M2 (blue curve) by the full range of frequency available for the GEM-2. The horizontal dashed line at 0.02 indicates the transition from low to moderate induction numbers [80]. Dots indicate the induction numbers for a set of six selectable operating frequencies among those currently available for the GEM-2 (minimum frequency = 30 Hz; maximum frequency = 93 kHz).

Regarding the in-phase component, although not visible in Figure 11 for drawing scale reasons, it should be noted that at low frequencies it is negative in all cases. For both models M1 and M2, in detail, the curves tend asymptotically towards values of −440 ppm for the one calculated with the sensor at 0.2 m above the ground, and −210 ppm for that at 0.9 m. As pointed out by Huang and Fraser (2003) [82] and by Farquharson et al. (2003) [83], as the frequency (or the induction number) takes small values, the complex response function becomes dominated by the magnetization effect, which is in-phase with and in the same direction as the primary magnetic field. This suggests that the negative values of the in-phase component at low frequencies, observed in Figure 11, are due to the susceptible materials present in the earth models (Figure 6a,b) or, more specifically, to the induced magnetization the susceptible materials exhibit when subjected to a magnetic field (no matter whether alternating or static). Here, we wish to highlight that the low-frequency asymptotic values depend exclusively on the magnetic permeability of the materials, unlike the values of the in-phase component at moderate and high frequencies, which, on the other hand, are influenced by electrical conductivity as well. This is a good reason to always include a very low frequency among those to be selected to set up a multi-frequency device. In addition, we want to emphasize that when the recorded in-phase data contain negative values, a careful direct modelling performed a posteriori may be particularly useful for data interpretation, using, in this case, an inversion algorithm taking into account both electrical conductivity and magnetic permeability simultaneously [83].

To quantify to what extent the complex EMI responses described above are affected by a modification in the value of the electrical conductivity and/or magnetic permeability of earth models M1 and M2, we have estimated a whole set of sensitivity functions for each device, using Equations (13) and (14) and assuming an operating height of 0.9 m. Figures 13 and 14 show the sensitivity functions as functions of electrical conductivity and magnetic permeability for both the quadrature (Q) and in-phase (P) components of the Dualem21H device above the two earth models M1 (Figures 13a–d and 14a–d) and M2 (Figures 13e–h and 14e–h). Similarly, Figures 15 and 16 show the sensitivities for the Explorer. Finally, Figure 17 shows, frequency by frequency, the sensitivities estimated for the GEM-2 above earth models M1 and M2. Note that in all graphs, the sensitivities are plotted in a non-normalized form, with the values expressed using the appropriate units of measurement for both numerator (Q or P component of the response, in mS/m or ppm for the former and in ppt or ppm for the latter) and denominator (electrical conductivity, in S/m, or magnetic permeability, in H/m) of Equations (13) and (14). We adopted this graphical representation because it allows users to quantitatively compare the whole set of sensitivity functions. It is the standard representation used in the FDEMtools3 package; however, users can modify some scripts to get other representations, similar to the ones shown in the Supplementary material (Figures S5–S8) in [7].

The analysis and comparison of the sensitivity functions in Figures 13–17 certainly provide further useful information to select the most appropriate and best configured measuring device to better characterize the target in the M1 and M2 models. Here, we leave this choice to the reader, according to their own analyses, comparisons, and considerations.

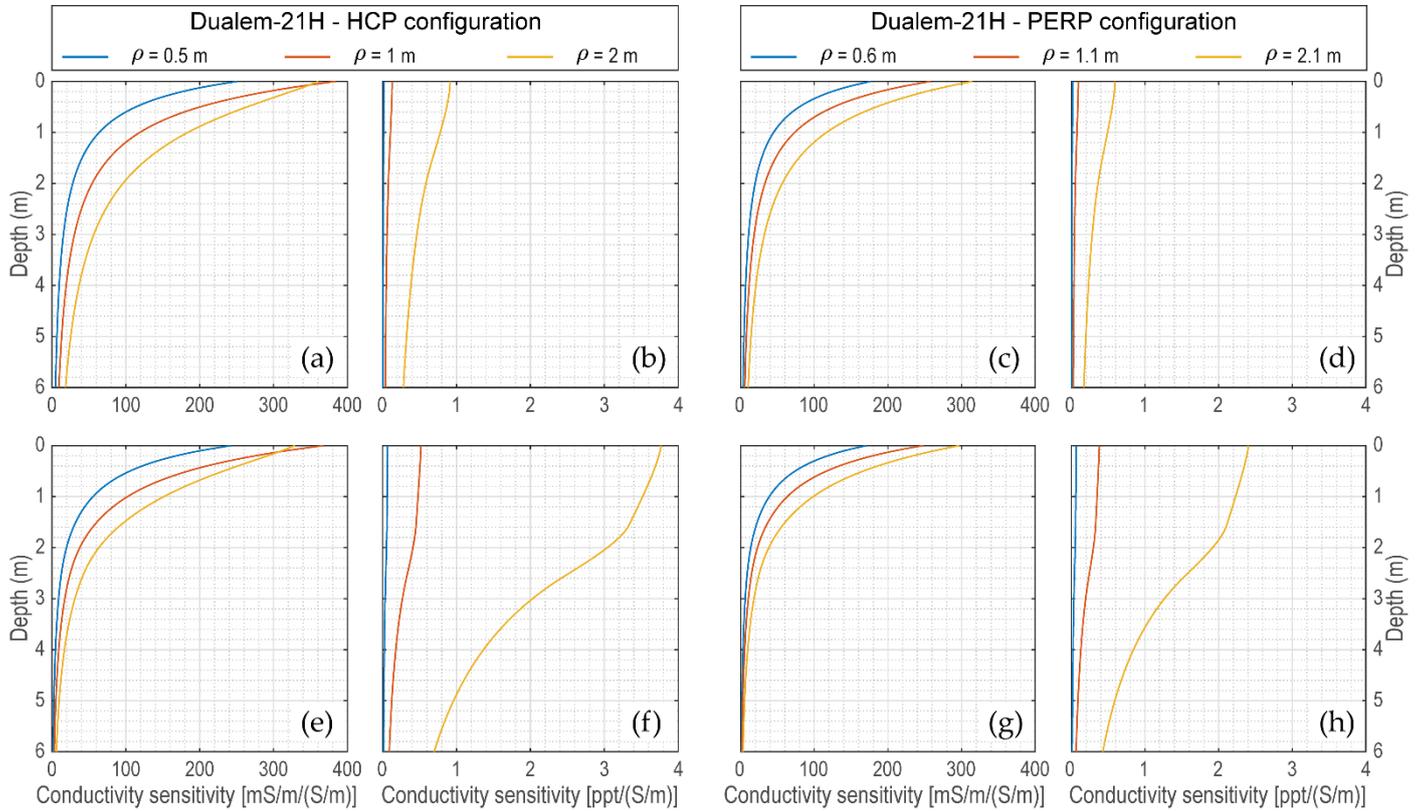

**Figure 13.** Sensitivity functions to electrical conductivity of the Dualem-21H. (**a,c**) Q and (**b,d**) P sensitivities at 0.9 m above model M1; (**e,g**) Q and (**f,h**) P sensitivities at 0.9 m above model M2.

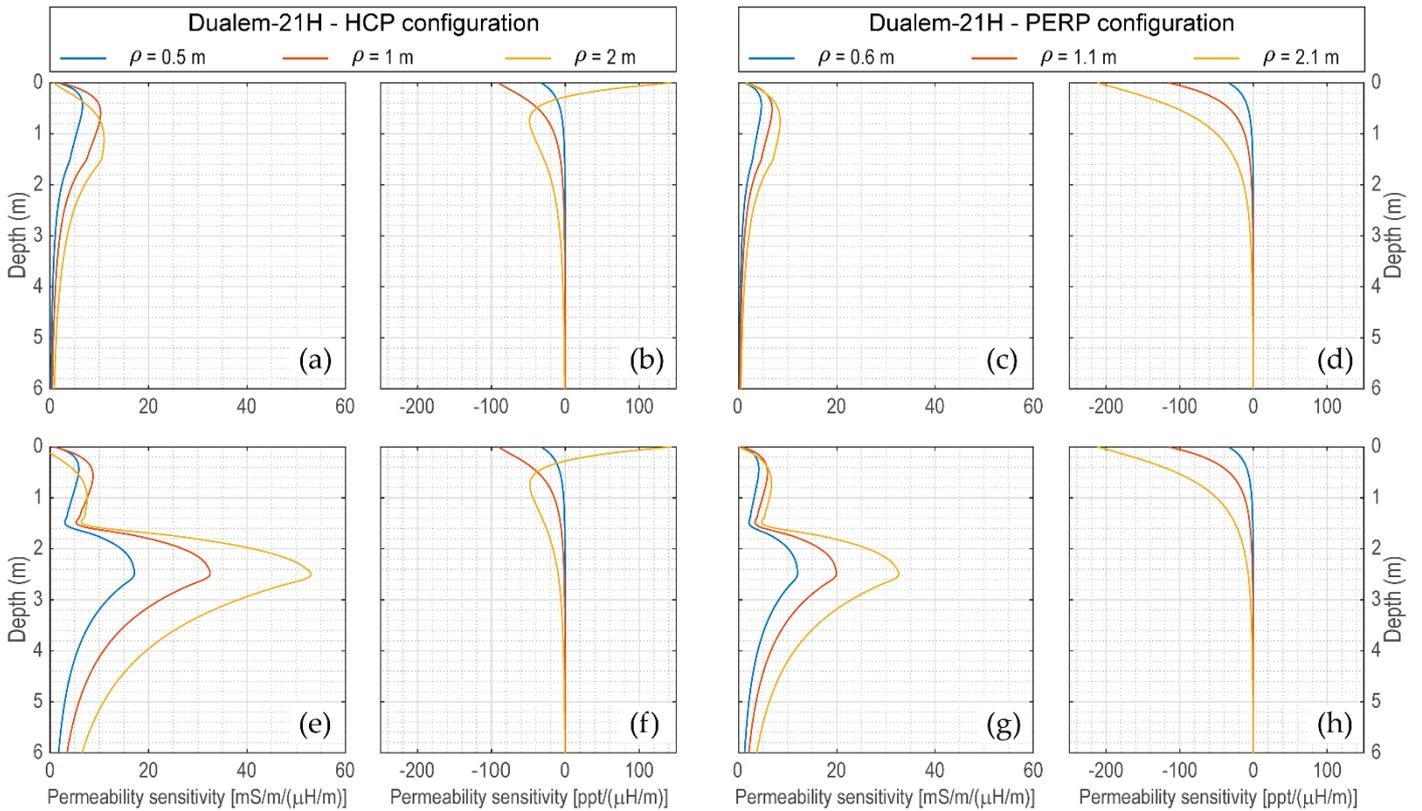

**Figure 14.** Sensitivity functions to magnetic permeability of the Dualem-21H. (**a**,**c**) Q and (**b**,**d**) P sensitivities at 0.9 m above model M1; (**e**,**g**) Q and (**f**,**h**) P sensitivities at 0.9 m above model M2.

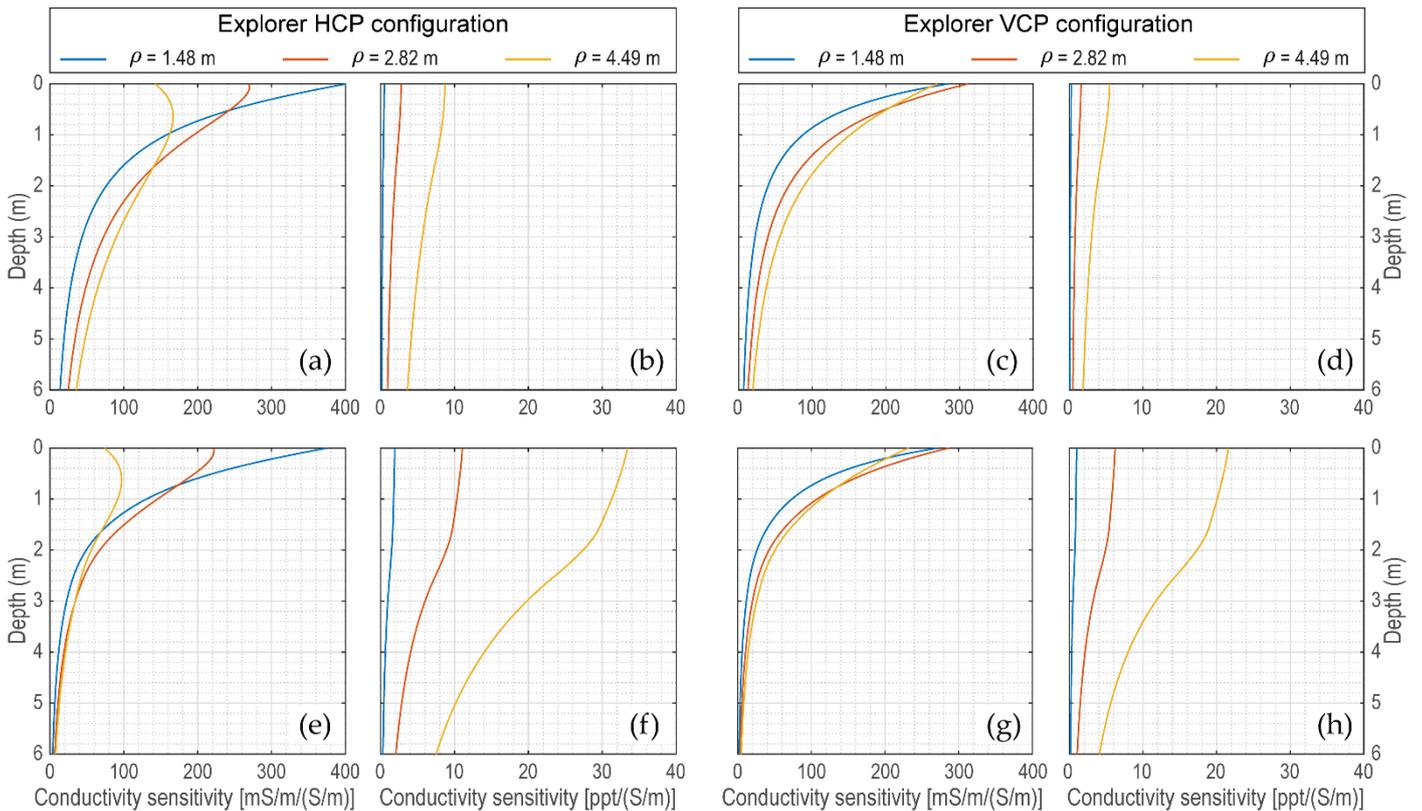

**Figure 15.** Sensitivity functions to electrical conductivity of the Explorer. (**a**,**c**) Q and (**b**,**d**) P sensitivities at 0.9 m above model M1; (**e**,**g**) Q and (**f**,**h**) P sensitivities at 0.9 m above model M2.

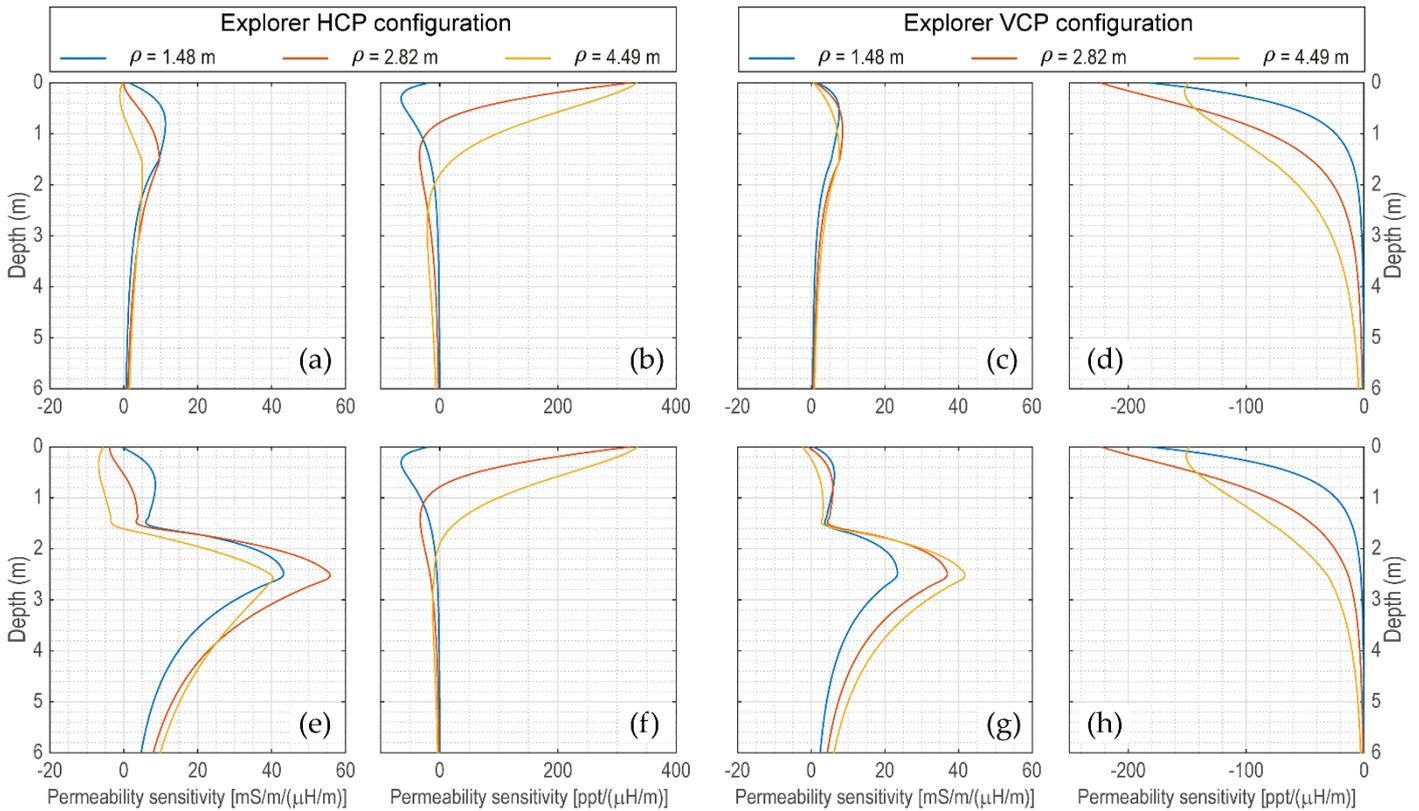

**Figure 16.** Sensitivity functions to magnetic permeability of the Explorer. (**a**,**c**) Q and (**b**,**d**) P sensitivities at 0.9 m above model M1; (**e**,**g**) Q and (**f**,**h**) P sensitivities at 0.9 m above model M2.

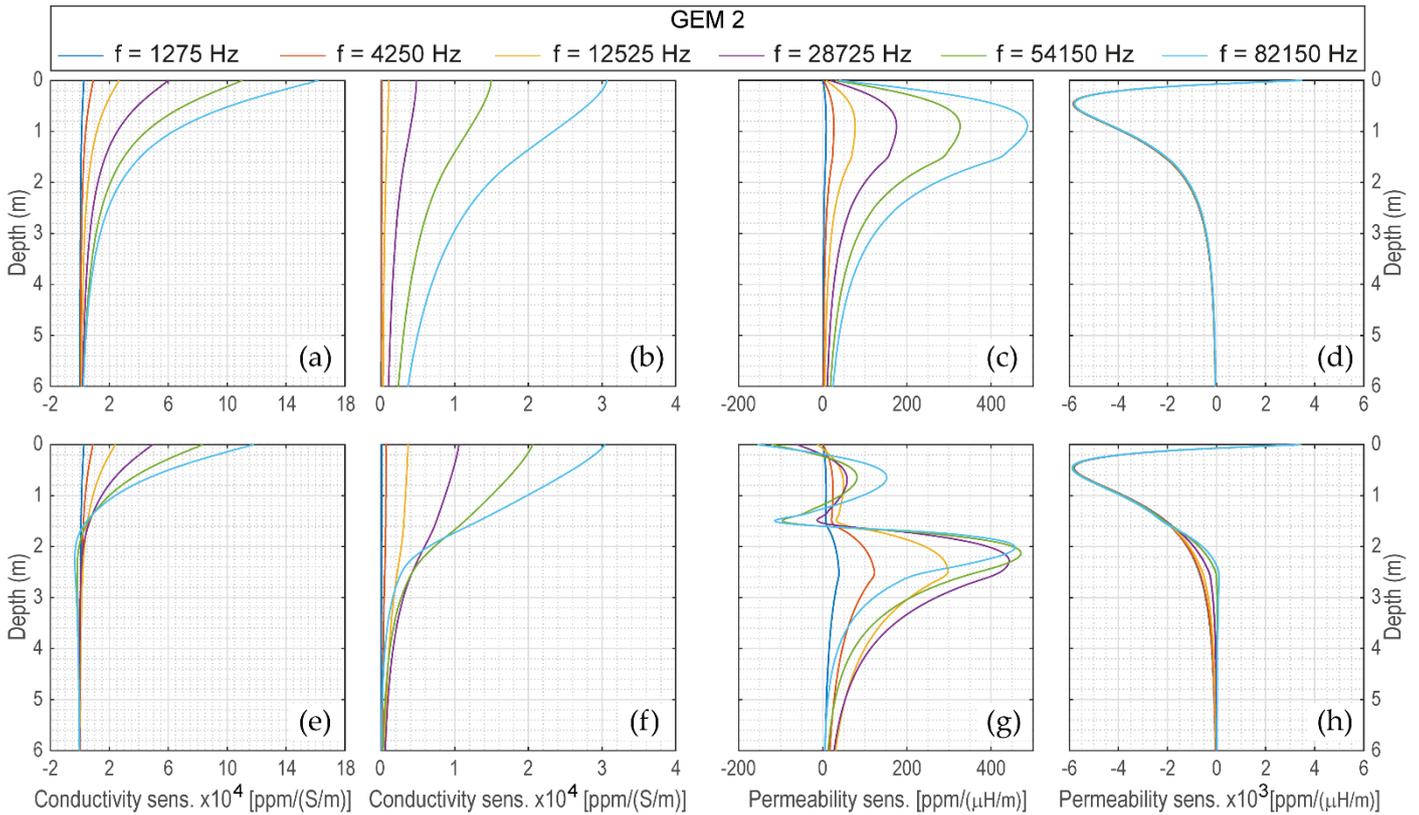

**Figure 17.** Sensitivity curves of the GEM-2. (**a**,**e**) Q and (**b**,**f**) P sensitivities to electrical conductivity at 0.9 m above models M1 and M2, respectively; (**c**,**g**) Q and (**d**,**h**) P sensitivities to magnetic permeability at 0.9 m above models M1 and M2, respectively.

Finally, to complete this numerical example, we report the values that one of the Auxiliary Routines (fdemdoi.m) of FDEMtools3 provides for the DOI (Table 4). As mentioned above (Section 2.4), these values are only indicative and useful to obtain an approximate estimate of the DOI, according to the criterion adopted in [67]. For each of earth models M1 or M2, Table 4 lists the DOIs achievable by each of the three devices, assuming an operating height of 0.9 m above the ground. Such values can also be graphically estimated by plotting the cumulative response functions, which are optional outputs of the forward modelling package in the FDEMtools3. Running the FDEMforward GUI with the "cumulative response graph" option activated, the cumulative response functions are firstly computed by integrating the sensitivity functions, and then plotted down to the depth that coincides with the DOI.

**Table 4.** Depth of Investigations (m) estimated with an operating height of 0.9 m.

| Model | Dualem-21H | | | | | | CMD Explorer | | | | | | GEM-2 | | | | | |
|---|---|---|---|---|---|---|---|---|---|---|---|---|---|---|---|---|---|---|
| | $^{HCP}\rho_1$ | $^{HCP}\rho_2$ | $^{HCP}\rho_3$ | $^{PERP}\rho_1$ | $^{PERP}\rho_2$ | $^{PERP}\rho_3$ | $^{HCP}\rho_1$ | $^{HCP}\rho_2$ | $^{HCP}\rho_3$ | $^{VCP}\rho_1$ | $^{VCP}\rho_2$ | $^{VCP}\rho_3$ | $f_1$ | $f_2$ | $f_3$ | $f_4$ | $f_5$ | $f_6$ |
| M1 | 2.7 | 3.5 | 6.7 | 2.7 | 2.8 | 3.9 | 5.5 | 5.9 | 7.9 | 3.2 | 5.1 | 8.3 | 8.3 | 8.3 | 8.3 | 8.3 | 8.3 | 7.9 |
| M2 | 2.7 | 3.5 | 6.3 | 2.4 | 2.8 | 3.9 | 5.1 | 5.6 | 7.1 | 3.2 | 4.7 | 7.5 | 8.3 | 8.3 | 7.2 | 6.4 | 5.5 | 4.7 |

$\rho_i$ (i = 1, 2, 3) are the inter-coil distances; $f_i$ (i = 1, 2, 3, 4, 5, 6) are the operating frequencies chosen for the example.

We remark here that all the Figures reported in this Section can be generated and displayed by running the test script drawfigures.m located in the main directory FDEMtools3.

## 5. Conclusions

A simulation of the model response to a prescribed distribution of the electromagnetic features in a stratified subsoil is a very useful tool to plan a data acquisition campaign and adopt an effective sensing device, with the best possible configuration. This is possible whenever a geophysicist has some a priori information on the surveying site and knows which physical target he/she is going to observe.

After recalling some basic concepts about the earth propagation of an electromagnetic field and discussing some physical quantities which are critical for the correct comprehension of the phenomenon, an interactive software tool for forward modelling has been introduced. It reproduces the model response to a given electromagnetic features distribution and allows the simulation of data acquisition by a specific instrument, either existing or hypothetical.

To illustrate the use of the package, two three-layer earth models have been analyzed by comparing the response of three commercial devices. The simulation shows that an effective choice of a specific sensing device, as well as its correct configuration, can only be performed by taking into consideration the target characteristics and the operating height of the device. At the same time, drawing sensitivity functions and cumulative response graphs is crucial to distinguish between data-driven and model-driven inversion results, and determine a reliable depth of investigation. A forward modelling software simulator turns out to be a precious tool to assist a geophysicist in planning a surveying session.

Besides expanding the FDEMtools3 toolbox by implementing a graphical user interface for forward modelling and the corresponding computational routines, this paper also introduces in the package a new regularized minimal-norm inversion algorithm, which helps in selecting a suitably regular solution for the underdetermined least-squares problems to be solved, and allows the use of a model profile for the solution, in those cases where such information is available.

**Author Contributions::** Conceptualization, G.P.D., P.D.d.A., F.P. and G.R.; methodology, G.P.D., P.D.d.A., F.P. and G.R.; software, G.P.D., P.D.d.A., F.P. and G.R.; validation, G.P.D., P.D.d.A., F.P. and G.R.; formal analysis, G.P.D., P.D.d.A., F.P. and G.R.; writing—original draft preparation, G.P.D., P.D.d.A., F.P. and G.R.; writing—review and editing, G.P.D., P.D.d.A., F.P. and G.R. All authors have read and agreed to the published version of the manuscript.

**Funding:** This research was partially funded by Fondazione di Sardegna, Progetto biennale bando 2021, "Computational Methods and Networks in Civil Engineering (COMANCHE)". F.P., P.D.A. and G.R. were partially supported by the INdAM-GNCS 2022 project "Metodi e modelli di regolarizzazione per problemi malposti di grandi dimensioni". P.D.A.

gratefully acknowledges Fondo Sociale Europeo REACT EU—Programma Operativo Nazionale Ricerca e Innovazione 2014–2020 and Ministero dell'Università e della Ricerca for the financial support.

## Appendix A: Brief review of the Maxwell equations

Electromagnetic induction phenomena obey Maxwell's equations, which describe how electric and magnetic fields are generated by charges, currents, and changes of the fields. The differential form in the time domain of these equations is given by

$$\nabla \cdot \boldsymbol{D} = q, \qquad \text{Gauss' law} \tag{A1}$$

$$\nabla \cdot \boldsymbol{B} = 0, \qquad \text{Gauss' law for magnetic fields} \tag{A2}$$

$$\nabla \times \boldsymbol{E} = -\frac{\partial \boldsymbol{B}}{\partial t}, \qquad \text{Faraday's law} \tag{A3}$$

$$\nabla \times \boldsymbol{H} = \boldsymbol{J} + \frac{\partial \boldsymbol{D}}{\partial t}, \qquad \text{Ampère-Maxwell's law} \tag{A4}$$

where $\boldsymbol{D}$ is the dielectric displacement (C/m$^2$), $\boldsymbol{B}$ the magnetic flux density or the magnetic induction (T), $\boldsymbol{E}$ the electric field intensity (V/m), $\boldsymbol{H}$ the magnetic field intensity (A/m), $\boldsymbol{J}$ the electric current density (A/m$^2$), and $q$ the electric charge density (C/m$^3$). The symbols $\nabla \cdot$ and $\nabla \times$ stand for divergence and curl operators, respectively. These equations are usually coupled through the following constitutive relations:

$$\boldsymbol{D} = \varepsilon \boldsymbol{E}, \tag{A5}$$

$$\boldsymbol{B} = \mu \boldsymbol{H}, \tag{A6}$$

$$\boldsymbol{J} = \sigma \boldsymbol{E}, \tag{A7}$$

where $\varepsilon$, $\mu$, $\sigma$, are the dielectric permittivity (F/m), the magnetic permeability (H/m), and the electric conductivity (S/m) of a conductive magnetic material. In free space, where the electric conductivity is zero, the dielectric permittivity and the magnetic permeability take the values $\varepsilon_0 = 8.854 \cdot 10^{-12}$ F/m and $\mu_0 = 4\pi \cdot 10^{-7}$ H/m, respectively. For any medium other than a vacuum, the ratio of the permeabilities of a medium to that of free space defines the dimensionless relative permeability $\mu_r = \frac{\mu}{\mu_0}$ as well the ratio $\varepsilon_r = \frac{\varepsilon}{\varepsilon_0}$ defines the relative dielectric permittivity.

For a magnetic material, Equation (A6) can be expressed in terms of a diagnostic parameter, the magnetic susceptibility $\chi$, which measures how much a material is susceptible to being magnetized.

In terms of relative permeability, it is

$$\chi = \mu_r - 1, \tag{A8}$$

so that the magnetic permeability is

$$\mu = \mu_0 + \mu_0 \chi. \tag{A9}$$

Therefore, it follows that the magnetic induction field, $\boldsymbol{B}$ (Equation (A6)), can be expressed as

$$\boldsymbol{B} = \mu \boldsymbol{H} = \mu_0 (1 + \chi) \boldsymbol{H} = \mu_0 \boldsymbol{H} + \mu_0 \chi \boldsymbol{H} = \mu_0 \boldsymbol{H} + \mu_0 \boldsymbol{M}, \tag{A10}$$

where

$$\boldsymbol{M} = \chi \boldsymbol{H} \tag{A11}$$

is the magnetization field that the material acquires when a magnetic field intensity, $\boldsymbol{H}$, acts on it.

As shown in [44], Maxwell's equations together with the constitutive relations can be combined to yield the electromagnetic wave equations for propagation (as wave and diffusion) of electric and magnetic fields in an isotropic homogeneous lossy medium having electric conductivity $\sigma$, magnetic permeability $\mu$, and dielectric permittivity $\varepsilon$. Taking the curl of Equation (A3) and using, in the following order, Equations (A6), (A4), (A7), and (A5), the use of the identity $\nabla \times \nabla \times \boldsymbol{E} = -\nabla^2 \boldsymbol{E}$, where the symbol $\nabla^2$ stands for the Laplacian operator, gives the equation for the electric field in time domain

$$\nabla^2 \boldsymbol{E} - \mu\sigma\frac{\partial \boldsymbol{E}}{\partial t} - \mu\varepsilon\frac{\partial^2 \boldsymbol{E}}{\partial t^2} = 0. \tag{A12}$$

Likewise, taking the curl of equation (A4) and using, in the following order, Equations (A7), (A3), (A5), and (A6), the use of the identity $\nabla \times \nabla \times \boldsymbol{H} = -\nabla^2 \boldsymbol{H}$ yields the equation for the magnetic field in time domain

$$\nabla^2 \boldsymbol{H} - \mu\sigma\frac{\partial \boldsymbol{H}}{\partial t} - \mu\varepsilon\frac{\partial^2 \boldsymbol{H}}{\partial t^2} = 0. \tag{A13}$$

Considering harmonically varying fields at angular frequency $\omega$, that is $\boldsymbol{E} = \boldsymbol{E}_0 e^{-i\omega t}$ and $\boldsymbol{H} = \boldsymbol{H}_0 e^{-i\omega t}$, Equations (A12) and (A13) become the two following Helmoltz equations

$$\nabla^2 \boldsymbol{E} + i\omega\mu\sigma\boldsymbol{E} + \omega^2\mu\varepsilon\boldsymbol{E} = \nabla^2 \boldsymbol{E} + k^2\boldsymbol{E} = 0, \tag{A14}$$

and

$$\nabla^2 \boldsymbol{H} + i\omega\mu\sigma\boldsymbol{H} + \omega^2\mu\varepsilon\boldsymbol{H} = \nabla^2 \boldsymbol{H} + k^2\boldsymbol{H} = 0, \tag{A15}$$

where

$$k = \sqrt{\omega^2\mu\varepsilon + i\omega\mu\sigma} = a + ib \tag{A16}$$

is the complex wavenumber, whose real and imaginary parts are respectively given by [46]:

$$a = \omega\sqrt{\frac{\mu\varepsilon}{2}\left(\sqrt{1 + \frac{\sigma^2}{\omega^2\varepsilon^2}} + 1\right)} \tag{A17}$$

and

$$b = \omega\sqrt{\frac{\mu\varepsilon}{2}\left(\sqrt{1 + \frac{\sigma^2}{\omega^2\varepsilon^2}} - 1\right)}. \tag{A18}$$

The imaginary part, which is also called the attenuation coefficient, plays a key role in electromagnetism since its inverse defines the skin depth $\delta$.

### Appendix A.1. Quasi-Stationary Approximation

Alternating electromagnetic fields that vary slowly with time are referred to as low-frequency alternating fields or quasi-stationary fields. In the case of quasi-stationary fields, Maxwell's equations can be simplified by dropping the term $\frac{\partial \boldsymbol{D}}{\partial t}$ in Ampère-Maxwell's law (Equation (A4)) but retaining the term $\frac{\partial \boldsymbol{B}}{\partial t}$ in Faraday's law (Equation (A3)). This means that the displacement current is negligible with respect to the conduction current, which remains the only source of the quasi-stationary magnetic field. This also means that the electromagnetic properties of the medium are such that $\sigma \gg \omega\varepsilon$. Then, Equations (A14) and (A15) can be approximated as

$$\nabla^2 \boldsymbol{E} + i\omega\mu\sigma\boldsymbol{E} \simeq 0 \tag{A19}$$

and

$$\nabla^2 \boldsymbol{H} + i\omega\mu\sigma\boldsymbol{H} \simeq 0, \tag{A20}$$

which are known as the diffusion equations of electromagnetic fields. They describe the penetration of electromagnetic fields (but do not consider wave propagation) in an isotropic homogeneous lossy medium having electric conductivity $\sigma$ and magnetic permeability $\mu$. The complex wavenumber $k$ becomes

$$k = \sqrt{i\omega\mu\sigma} = a + ib, \tag{A21}$$

whose real and imaginary parts are

$$a = b = \sqrt{|k^2|} = \sqrt{\frac{\omega\mu\sigma}{2}} = \frac{1}{\delta}, \tag{A22}$$

where

$$\delta = \sqrt{\frac{2}{\omega\mu\sigma}} \tag{A23}$$

is the skin depth.

## Appendix B: Step-by-step electromagnetic induction

### Appendix B.1. Step 1

Let us first consider two nearby coils in free space (or in free air), as in Figure A1a. Suppose that coil Tx (Transmitter) is connected to an external alternating voltage source, while coil Rx (Receiver) is connected to a voltmeter to read voltages in it (Figure A1b). Let

$$I_P = I_0 \cdot e^{i\omega t} \tag{A24}$$

the sinusoidal alternating current driven in the primary coil by the external voltage source. According to Ampère-Maxwell's law, this current produces a time-varying magnetic field intensity, $\mathbf{H}_P$, or magnetic flux density, $\boldsymbol{B}_P = \mu_0 \boldsymbol{H}_P$, around the loop, which alternates with the same frequency and phase as the current (Figure A1c). Both magnitude and direction of this field vary with position in a complex way around the coil, but its magnitude is always proportional to the current flowing in the coil: $|\boldsymbol{B}_P| \propto I_P$.

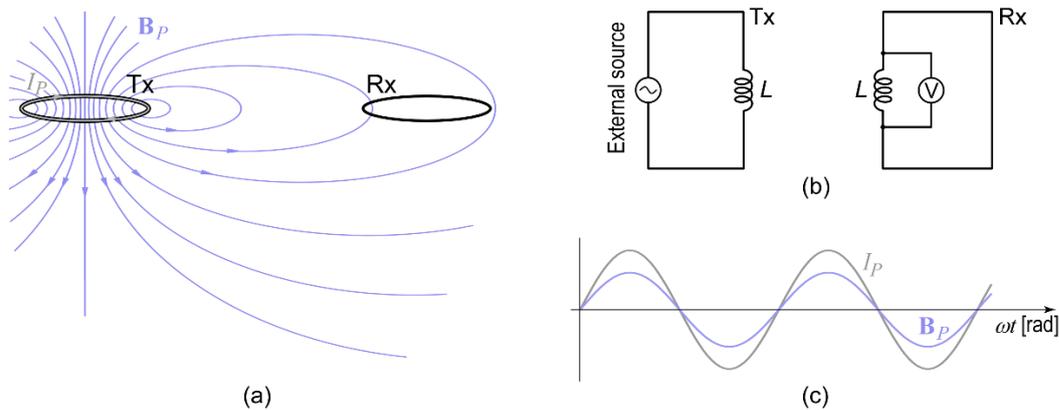

**Figure A1.** (**a**) Sketch of two magnetically coupled coils in a free space. (**b**) Equivalent single-loops circuits for the transmitter (**on the left**) and the receiver (**on the right**) coils. (**c**) Primary current and primary magnetic field as a function of time.

The time-varying magnetic field generates a changing magnetic flux through coil Rx, $\Phi_R(\boldsymbol{B}_P)$. Therefore, the magnetic field interacts with coil Rx to produce an electromotive force, according to Faraday's law:

$$\mathcal{E} = -\frac{\partial \Phi_R(\boldsymbol{B}_P)}{\partial t}. \tag{A25}$$

Since the magnetic field is proportional to the current $I_P$, and the magnetic flux, by definition, is proportional to the magnetic field, the magnetic flux through coil Rx is proportional to the current flowing in coil Tx, that is,

$$\Phi_R(\boldsymbol{B}_P) = M_{TR} I_P, \tag{A26}$$

where $M_{TR}$ is the mutual inductance, which is defined as the magnetic flux that passes through coil Rx due to a unit electric current circulating in coil Tx. The mutual inductance $M_{TR}$ depends on the geometry of the coils, their relative orientation and distance, and on the magnetic permeability of free space $\mu_0$ ($4\pi \cdot 10^{-7}$ H/m). Combining Equations (A25) and (A26), the voltage sensed by coil Rx is

$$\mathcal{E}_{TR} = -\frac{\partial \Phi_R(\boldsymbol{B}_P)}{\partial t} = -M_{TR}\frac{\partial I_P}{\partial t} = -i\omega M_{TR} I_P. \tag{A27}$$

This voltage is usually employed to measure the primary magnetic field at the receiver.

## *Appendix B.2. Step 2*

Now, let us consider again the two coils Tx and Rx in free air but above a half-space containing a conductive magnetic body with electrical conductivity $\sigma$ and magnetic permeability $\mu$ (Figure A2a).

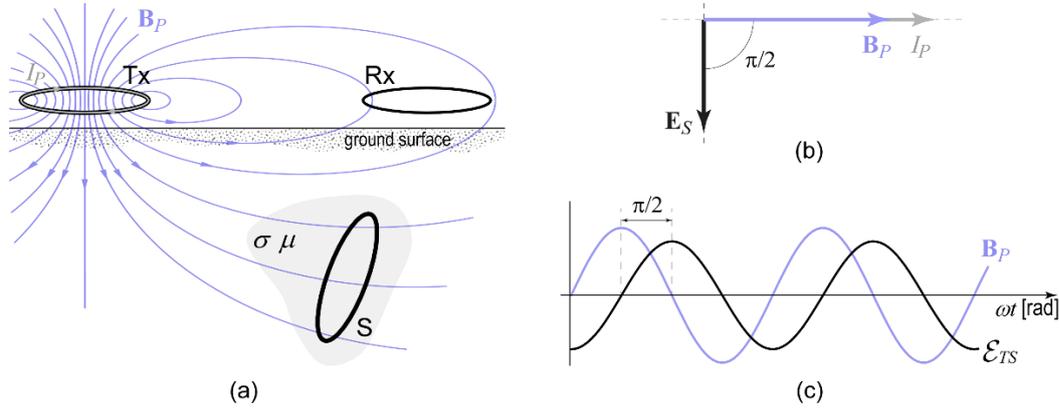

**Figure A2.** (**a**) Sketch of two magnetically coupled coils over a half-space with a conductive body. (**b**) Phasor diagram using as a reference the primary current. The secondary electric field lags the primary magnetic field by 90°. (**c**) Time dependence of $\mathscr{E}_{TS}(t)$ and $\boldsymbol{B}_P(t)$ that shows the same lagging phase as in the phasor diagram.

For a bulk material (the conductive magnetic body) there is not a loop per se, but many short-circuited loops. However, Faraday's law is general and it does not require the existence of a physical loop. Faraday's law states that when the magnetic flux through a surface changes, a time-varying electric field is induced along the boundary of that surface. This is true for any closed loop, either in empty space or in a physical material, through which the magnetic flux is changing over time. Thus, assuming S as one of these loops inside the body (Figure A2a), the standard integral form of Faraday's law reads

$$\oint_S \boldsymbol{E}_S \cdot d\boldsymbol{l} = -\frac{\partial \Phi(\boldsymbol{B_P})}{\partial t}, \tag{A28}$$

where $\mathbf{E}_S$ is the electric field at every point of such a loop and $d\mathbf{l}$ is an oriented displacement along the loop. The induced electromotive force $\mathscr{E}_{TS}$ is related to $\mathbf{E}_S$ by

$$\mathscr{E}_{TS} = \oint_S \boldsymbol{E}_S \cdot d\boldsymbol{l}. \tag{A29}$$

Therefore, as in the case of coil Rx (Equation (A27), by introducing the mutual induction $M_{TS}$ the electromotive force induced in the loop S can be expressed in terms of the primary current $I_P$ by

$$\mathscr{E}_{TS} = -\frac{\partial \Phi(\boldsymbol{B_P})}{\partial t} = -M_{TS}\frac{\partial I_P}{\partial t} = -i\omega M_{TS} I_P. \tag{A30}$$

This electromotive force alternates with the same frequency as the primary current but lagging by 90° behind the current (or the primary magnetic field); see Figure A2c. The mutual inductance $M_{TS}$ depends on the geometry of coils Tx and S, on their relative orientation and distance, and on the magnetic permeability $\mu$ of the core material in loop S.

## *Appendix B.3. Step 3*

The alternating voltage induced in the conductive body by the time-varying primary magnetic field causes alternating currents to flow in the bulk material as they do through wires. These are the eddy currents that flow along closed loops concentrated near the boundary surface of the body (skin effect) and in planes perpendicular to the magnetic field causing them. Let S be one of these closed loops (Figure A3a). Figure A3b shows its equivalent single-loop circuit with lumped resistance $R$ and inductance $L$. Let $\mathscr{E}_{TS}(t)$ be the alternating voltage source that establishes the alternating current, $I_{eddy}$.

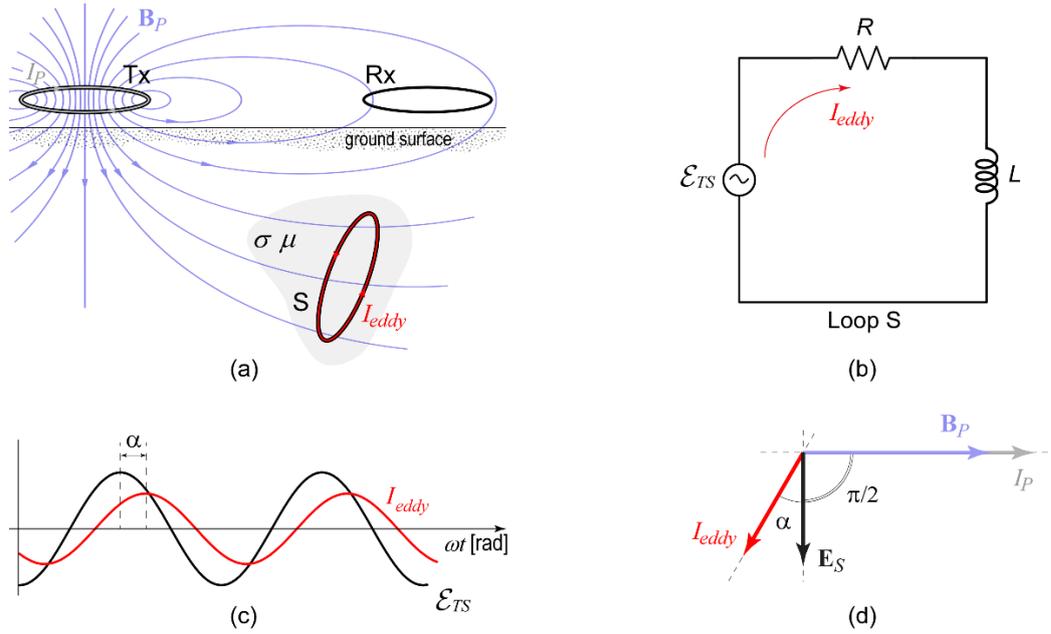

**Figure A3.** (**a**) Sketch showing one eddy current loop. (**b**) Equivalent circuit with lumped resistance R and inductance L, connected across an alternating voltage source. (**c**) Time dependence of $\mathcal{E}_{TS}(t)$ and $I_{eddy}(t)$ across loop S; the current lags the voltage. (**d**) Phasor diagram using as a reference the primary current. The secondary electric field lags by 90° the primary magnetic field, while eddy currents show an additional phase lag.

Applying Kirchhoff's voltage rule, the circuit equation reads

$$\mathcal{E}_{TS} - R I_{eddy} - L \frac{dI_{eddy}}{dt} = 0, \tag{A31}$$

which, for the present time-harmonic case, yields

$$\mathcal{E}_{TS} = (R + i\omega L) I_{eddy}. \tag{A32}$$

The complex quantity in the brackets is the impedance of the $RL$ circuit, whose amplitude is given by

$$|Z| = \sqrt{R^2 + \omega^2 L^2}, \tag{A33}$$

which, for the present time-harmonic case, yields the phase

$$\alpha = arctan\left(\frac{\omega L}{R}\right). \tag{A34}$$

Therefore, letting $\mathcal{E}_{TS}(t) = \mathcal{E}_0 \cdot e^{i\left(\omega t - \frac{\pi}{2}\right)}$, the sinusoidal alternating current circulating in the circuit is

$$I_{eddy}(t) = \frac{\mathcal{E}_0}{\sqrt{R^2 + \omega^2 L^2}} \cdot e^{i(\omega t - \frac{\pi}{2} - \alpha)}, \tag{A35}$$

which lags the voltage by $\alpha$ radians (Figure A3c) and the primary magnetic field (or primary current) by $\alpha + \frac{\pi}{2}$ radians (Figure A3d).

The phase shift $\alpha$ depends only on the response parameter $\beta = \omega \frac{L}{R}$, also known as the dimensionless induction number. When $\beta \to 0$ or equivalently $R \to \infty$ (for a given value of $\omega L$), the circuit becomes purely resistive as the amplitude and phase of the impedance becomes $|Z| = R$ and $\alpha = 0$, respectively. In this case, the current circulating in the circuit is in-phase with the induced voltage $\mathcal{E}_{TS}(t)$ and is given by

$$I_{eddy}(t) = \frac{\mathcal{E}_0}{R} \cdot e^{i\left(\omega t - \frac{\pi}{2}\right)}. \tag{A36}$$

When $\beta \to \infty$ or equivalently $R \to 0$ (for a given value of $\omega L$), the circuit becomes purely inductive as the amplitude impedance takes the value $|Z| = \omega L$ and the phase approaches $\frac{\pi}{2}$ radians:

$$\alpha = \lim_{R \to 0}\left[arctan\left(\frac{\omega L}{R}\right)\right] = \frac{\pi}{2}. \tag{A37}$$

In this case, thus, the current circulating in the circuit is in quadrature with the induced voltage $\mathcal{E}_{TS}(t)$, lags the primary current by $\pi$ radians, and is given by

$$I_{eddy}(t) = \frac{\mathcal{E}_0}{\omega L} \cdot e^{i(\omega t - \pi)}. \tag{A38}$$

*Appendix B.4. Step 4*

Eddy currents induced in the body generate a time-varying magnetic field around the body (Figure A4a), according to Ampère-Maxwell's law and (A32). This field, called secondary magnetic field, generates in turn a secondary voltage in coil Rx, according to Faraday's law:

$$\mathscr{E}_{SR} = -i\omega M_{SR} I_{eddy} = -i\omega M_{SR}\frac{\mathscr{E}_{SR}}{R+i\omega L} = -\frac{\omega^2 M_{TS}M_{SR}}{R+i\omega L} \cdot I_P.$$ (A39)

where $M_{SR}$ denotes the mutual inductance between the coils S and Rx.

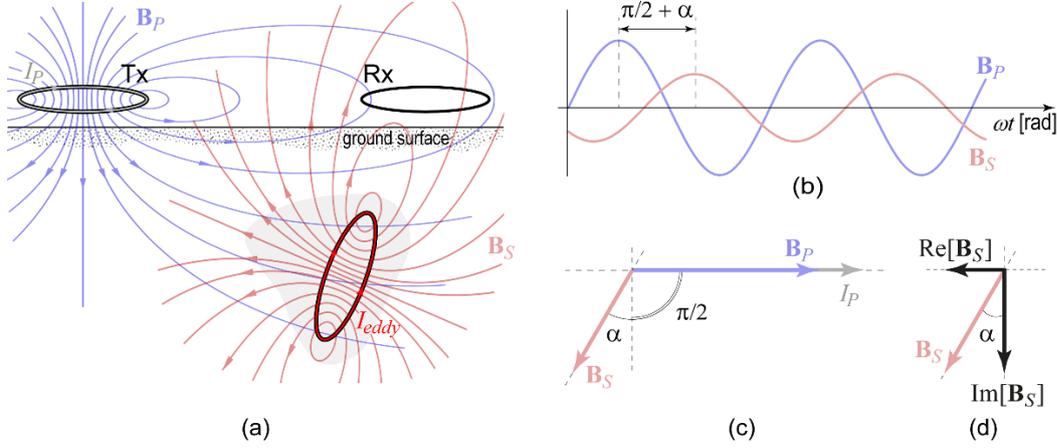

**Figure A4.** (**a**) Sketch showing the secondary magnetic field, $\mathbf{H}_S$, generated by the eddy current loop. (**b**) Time dependence of $\mathbf{H}_P$ and $\mathbf{H}_S$; the secondary magnetic field lags the primary field. (**c**) Phasor diagram using as a reference the primary current: the secondary magnetic field lags the primary magnetic field by $\alpha + \frac{\pi}{2}$ radians. (**d**) Decomposition of the secondary magnetic field in its real and imaginary parts, which are the in-phase and in quadrature components with respect to the primary field, respectively.

The receiver, then, simultaneously senses both the primary and the secondary magnetic fields, measuring both primary and secondary electromotive forces. In particular, the receiver records the whole electromagnetic response of the buried loop as the ratio of the secondary to the primary magnetic fields, which is equal to the ratio of the secondary to the primary voltages:

$$\frac{\mathscr{E}_S}{\mathscr{E}_P} = -\frac{M_{TS}\cdot M_{SR}}{M_{PR}\cdot L}\cdot\frac{i\beta}{1+i\beta} = \kappa\cdot\frac{i\beta}{1+i\beta} = \kappa\left(\frac{\beta^2+i\beta}{1+\beta^2}\right).$$ (A40)

The first factor

$$\kappa = -\frac{M_{TS}\cdot M_{SR}}{M_{PR}\cdot L},$$ (A41)

is the coupling coefficient. It depends only on relative size, shape, position, and orientation of the coils. The other factor, called the response function, is a complex-valued function of $\beta$, which depends on the frequency $\omega$ and on the target's electromagnetic properties:

$$G(\beta) = \frac{i\beta}{1+i\beta} = \frac{\beta^2}{1+\beta^2} + i\frac{\beta}{1+\beta^2}.$$ (A42)

Therefore, the electromagnetic response of the measuring device to the buried body is given by

$$M = \frac{\mathscr{E}_S}{\mathscr{E}_P} = \kappa\cdot G(\beta),$$ (A43)

$$Re\,M = \kappa\cdot\frac{\beta^2}{1+\beta^2},$$ (A44)

$$Im\,M = \kappa\cdot\frac{\beta}{1+\beta^2}.$$ (A45)

The real part of the response, having the same phase as the primary magnetic field, is usually designated the In-phase component, while the imaginary part, or Quadrature com-ponent, is out-of-phase with the primary by 90° (Figure A4c).

The response function becomes purely real when $\beta \to \infty$ (inductive limit), and when the instrument works at high frequency, or the target is highly conductive (low R) or highly inductive. Otherwise, the response function is purely imaginary when $\beta \to 0$ (resistive limit), which means using a low frequency, or being in the presence of a poorly conductive target (high R). Figure A5 shows the graph of the real and imaginary parts of the response function $G(\beta)$.

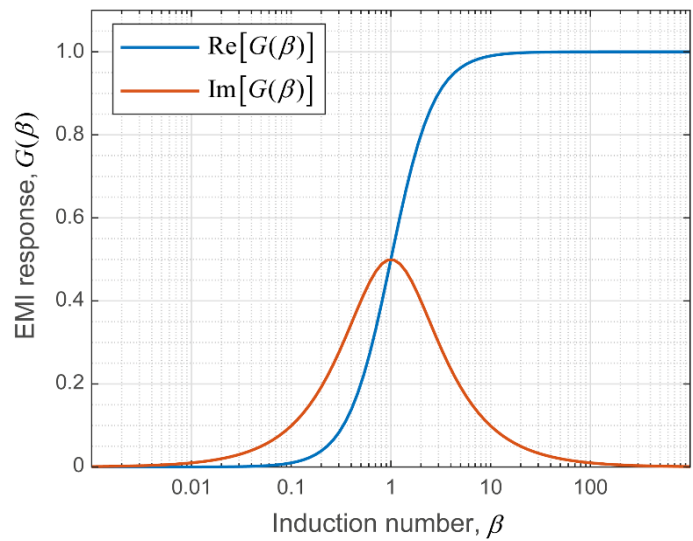

**Figure A5.** Real and Imaginary parts of the electromagnetic response function.